\documentclass[11pt]{amsart}
\usepackage{amssymb,amsmath,amsthm}
\usepackage{amscd}
\input amssym.def
\input amssym
\newtheorem{theorem}{Theorem}[section]

\newtheorem{lemma}{Lemma}[section]

\newtheorem{remark}{Remark}[section]
\newtheorem{definition}{Definition}[section]
\newtheorem{proposition}{Proposition}[section]

\newcommand{\ki}{\Lambda_{k,i}}

\newcommand{\be}{\begin{equation}}
\newcommand{\ee}{\end{equation}}

\renewcommand{\theequation}{\thesection.\arabic{equation}}
\renewcommand{\thetheorem}{\thesection.\arabic{theorem}}
\renewcommand{\theequation}{\thesection.\arabic{equation}}
\begin{document}

\title[Vertex-algebraic structure of principal subspaces, II]
{Vertex-algebraic structure of the principal subspaces of certain
$A_1^{(1)}$-modules, II: higher level case}

\author{C. Calinescu, J. Lepowsky and A. Milas}

\thanks{C.C. gratefully acknowledges partial support {}from the Center for
Discrete Mathematics and Theoretical Computer Science (DIMACS),
Rutgers University. J.L.  gratefully acknowledges partial support
{}from NSF grant DMS--0401302.}

\begin{abstract}
We give an a priori proof of the known presentations of (that is,
completeness of families of relations for) the principal subspaces of
all the standard $A_1^{(1)}$-modules.  These presentations had been
used by Capparelli, Lepowsky and Milas for the purpose of obtaining
the classical Rogers-Selberg recursions for the graded dimensions of
the principal subspaces.  This paper generalizes our previous paper.
\end{abstract}

\maketitle

\renewcommand{\theequation}{\thesection.\arabic{equation}}
\renewcommand{\thetheorem}{\thesection.\arabic{theorem}}
\setcounter{equation}{0} \setcounter{theorem}{0}
\setcounter{section}{0}

\section{Introduction}

The affine Kac-Moody algebra $A_1^{(1)}=\widehat{\goth{sl}(2)}$ is the
simplest infinite-dimensional Kac-Moody Lie algebra, and in some sense
the most prominent one.  Not only does $\widehat{\goth{sl}(2)}$ give
insight into the higher-rank affine Lie algebras, but in fact,
considerations of standard (= integrable highest weight)
$\widehat{\goth{sl}(2)}$-modules have frequently led to new ideas. For
instance, explicit constructions of the standard
$\widehat{\goth{sl}(2)}$-modules have been used to obtain
vertex-operator-theoretic derivations of the classical
Rogers-Ramanujan identities and related $q$-series identities
(cf. \cite{LW1}--\cite{LW4}, \cite{LP1}, \cite{LP2}, \cite{MP1},
\cite{MP2}). Another important use of standard
$\widehat{\goth{sl}(2)}$-modules is in the ``coset'' construction of
unitary Virasoro-algebra minimal models \cite{GKO}.  These
developments are deeply related to two-dimensional conformal field
theory.

More recently, to each standard $\widehat{\goth{sl}(n)}$-module
$L(\Lambda)$, Feigin and Stoyanovsky associated a distinguished
subspace $W(\Lambda)$, which they called the ``principal subspace'' of
$L(\Lambda)$ (\cite{FS1}, \cite{FS2}), and interestingly, the graded
dimensions of the principal subspaces of the standard
$\widehat{\goth{sl}(2)}$-modules are essentially the Gordon-Andrews
$q$-series (\cite{FS1}, \cite{Ge1}; cf. \cite{A}). These $q$-series
had previously appeared as the graded dimensions of the ``vacuum''
subspaces, with respect to a certain twisted Heisenberg subalgebra, of
the odd-level standard $\widehat{\goth{sl}(2)}$-modules
(\cite{LW2}--\cite{LW4}, \cite{MP1}). Since each standard
$\widehat{\goth{sl}(n)}$-module $L(\Lambda)$ of level $k$, $k \geq 1$,
is a module for a certain vertex operator algebra (\cite{FZ};
cf. \cite{DL}, \cite{LL}), it is natural to employ ideas {}from vertex
operator algebra theory to gain a better insight into the structure of
principal subspaces. In \cite{CLM1}--\cite{CLM2}, for the case
$\widehat{\goth{sl}(2)}$, the theory of vertex algebras and related
algebraic structures, including intertwining operators \cite{FHL}, has
been used to do this, via the construction of certain exact sequences,
which led to a vertex-algebra-theoretic interpretation of the
classical Rogers-Ramanujan and Rogers-Selberg recursions.  This in
turn explained the appearance of the Gordon-Andrews $q$-series, and
these $q$-series can be implemented by means of ``combinatorial
bases'' of the principal subspaces, revealing a fundamental
``difference-two condition'' that had already arisen in the setting of
\cite{LW2}--\cite{LW4}.

An important technical result used in \cite{CLM1} and \cite{CLM2} was
a certain presentation of (that is, the completeness of a certain
family of relations for) the principal subspaces of the standard
$\widehat{\goth{sl}(2)}$-modules (cf. Theorem 2.1 in \cite{CLM2}).
This result had been stated as Theorem 2.2.1$'$ in \cite{FS1}.
However, the proofs of this result that we are aware of all turn out
to require either a priori knowledge of a combinatorial basis of the
principal subspace $W(\Lambda)$ (see (\ref{ps}) below) or information
closely related to such knowledge.  But what one ideally wants is
rather an a priori proof of the presentation, which could then be used
to construct the exact sequences mentioned above, and thereby to
produce the bases. Thus it is an important problem to try to find an a
priori proof of the presentation of $W(\Lambda)$, and we were able to
achieve this for the level one standard
$\widehat{\goth{sl}(2)}$-modules in \cite{CalLM1}. Our proof in
\cite{CalLM1} was obtained in two steps. We first argued that the
presentation of $W(\Lambda_1)$ follows from the presentation of
$W(\Lambda_0)$, and then we proved the presentation of $W(\Lambda_0)$.
(These two steps are in fact interchangeable, so we could have placed
the proof of the presentation of $W(\Lambda_0)$ first.)

In the present paper we give an a priori proof of the presentation of
the principal subspaces more generally for all the standard
$\widehat{\goth{sl}(2)}$-modules.  The higher-level case brings
additional subtleties, and our approach is different from that in
\cite{CalLM1}.  Instead of trying to reduce the problem of proving the
presentation of principal subspaces to a ``preferred'' principal
subspace (e.g., $W(k \Lambda_0)$), we found it more convenient and
more elegant to prove the presentation of all the principal subspaces
of a given level at once.  This is done in the proofs of Theorems
\ref{th1} and \ref{th2} through a (necessarily) rather delicate
argument, which uses various properties of principal subspaces and
intertwining operators among standard modules. In our new approach all
the principal subspaces are on more-or-less equal footing.  Thus we
not only generalize the main result in \cite{CalLM1} to all the
standard $\widehat{\goth{sl}(2)}$-modules, but we also give a new
proof of the presentation of the principal subspaces in the level one
case, different from the one in \cite{CalLM1}.  This is why we write
the proof of the $k=1$ case separately and in full detail below; we
will also be generalizing this $k=1$ proof in a different direction
elsewhere.

This paper brings our in-depth analysis of the principal subspaces of
the standard $\widehat{\goth{sl}(2)}$-modules to an end.  Even though
the study of the principal subspaces of the
$\widehat{\goth{sl}(2)}$-modules is facilitated by the commutativity
of the underlying nilpotent Lie algebra used to define these
subspaces, many methods in this paper can be applied to more general
affine Lie algebras, both untwisted and twisted. In a sequel
\cite{CalLM2} we will shift our attention to standard modules for
affine Lie algebras of types $A,D,E$, in which case the relevant
nilpotent Lie algebras are nonabelian.

This paper is organized as follows. Section 2 gives the setting. In
Sections 3 and 4 we state and prove our main result, Theorem
\ref{th1}, which we also reformulate as Theorem \ref{th2}.  As in
\cite{CalLM1}, finding a further reformulation of the presentation of
the principal subspaces in terms of ideals of vertex (operator)
algebras is a natural problem.  This is achieved in Section 5 (Theorem
\ref{th3}), at least for the principal subspaces stemming {}from the
``vacuum'' higher-level $\widehat{\goth{sl}(2)}$-modules.

\section{The setting}

We start by recalling some background from \cite{CalLM1}, for the
reader's convenience.  Set
$$
\goth{g}=\goth{sl}(2)= \mathbb{C}x_{-\alpha} \oplus \mathbb{C} h
\oplus \mathbb{C} x_{\alpha},
$$
with bracket relations
$$
[h, x_{\alpha}]= 2 x_{\alpha}, \; \; [h, x_{-\alpha}]= -2
x_{-\alpha}, \; \; [x_{\alpha}, x_{-\alpha}]= h.
$$
The symmetric invariant bilinear form $\langle a, b \rangle= \mbox{tr}
\; (ab)$ ($a,b \in \goth{g}$) allows us to identify the Cartan
subalgebra $\goth{h}=\mathbb{C} h$ with its dual $\goth{h}^{*}$.  The
simple root $\alpha \in \goth{h}^{*}$ corresponding to the root vector
$x_{\alpha}$ identifies with $h \in \goth{h}$, that is, $h = \alpha$,
and $\langle \alpha, \alpha \rangle =2$.  Set $\goth{n}=\mathbb{C}
x_{\alpha}$.

We shall use the affine Lie algebra
\begin{equation}
\widehat{\goth{g}}= \goth{g} \otimes \mathbb{C}[t, t^{-1}] \oplus
\mathbb{C} {\bf k},
\end{equation}
with brackets
\begin{equation} \label{brackets}
[a \otimes t^m, b \otimes t^n]=[a,b] \otimes t^{m+n} + m \langle a, b
\rangle \delta_{m+n}{\bf k}
\end{equation}
for $a, b \in \goth{g}$, $m,n \in \mathbb{Z}$, with ${\bf k}$
central, and its subalgebras
\begin{eqnarray}
\bar{\goth{n}}&=&\mathbb{C}x_{\alpha} \otimes \mathbb{C}[t,
t^{-1}], \nonumber \\
\bar{\goth{n}}_{-}&=& \mathbb{C} x_{\alpha} \otimes
t^{-1}\mathbb{C}[t^{-1}], \nonumber \\
\bar{\goth{n}}_{\leq -2}& =&\mathbb{C} x_{\alpha} \otimes t^{-2}
\mathbb{C}[t^{-1}]. \nonumber \end{eqnarray}
The Lie algebra
$\widehat{\goth{g}}$ has the triangular decompositions
\begin{equation}
\widehat{\goth{g}}=( \mathbb{C} x_{-\alpha} \oplus \goth{g} \otimes
t^{-1} \mathbb{C}[t^{-1}]) \oplus (\goth{h} \oplus \mathbb{C}{\textbf k})
\oplus (\mathbb{C}x_{\alpha} \oplus \goth{g} \otimes t\mathbb{C}[t])
\end{equation}
and
\begin{equation}
\widehat{\goth{g}}= \widehat{\goth{g}}_{<0} \oplus
\widehat{\goth{g}}_{\geq 0},
\end{equation}
where
$$
\widehat{\goth{g}}_{<0}= \goth{g} \otimes t^{-1} \mathbb{C} [t^{-1}]
$$
and
$$ \widehat{\goth{g}}_{\geq 0}= \goth{g} \otimes \mathbb{C}[t] \oplus
\mathbb{C}{\textbf k}.
$$

Let $k \geq 1$.  We consider the level $k$ standard
$\widehat{\goth{g}}$-modules $L((k-i)\Lambda_0+i \Lambda_1)$,
where $\Lambda_0, \Lambda_1 \in (\goth{h} \oplus
\mathbb{C}{\textbf k})^*$ are the fundamental weights of
$\widehat{\goth{g}}$ ($\Lambda_j ({\textbf k}) =1$,
$\Lambda_j (h) =\delta_{j,1}$ for $j=0,1$) and $0
\leq i \leq k$ (cf. \cite{K}), so that $L(\Lambda_0)$ and
$L(\Lambda_1)$ are the level $1$ standard
$\widehat{\goth{g}}$-modules used in \cite{CalLM1}. For such
$i$, we set
\begin{equation}\label{Lambdakinotation}
\Lambda_{k,i}=(k-i) \Lambda_0+i \Lambda_1.
\end{equation}
Denote by $v_{\Lambda_{k,i}}$ a highest weight vector of
$L(\Lambda_{k,i})$.  (These highest weight vectors will be normalized
in Section 4 below.)

Throughout this paper we will write $x(m)$ for the action of $x
\otimes t^m \in \widehat{\goth{g}}$ on any
$\widehat{\goth{g}}$-module, where $x \in \goth{g}$ and $m \in
\mathbb{Z}$. In particular, we have the operator $x_{\alpha}(m)$, the
image of $x_{\alpha} \otimes t^m$. Sometimes we will simply write
$x(m)$ for the Lie algebra element $x \otimes t^m$. It will be clear
{}from the context whether $x(m)$ is an operator or a Lie algebra
element.

We generalize the definition of the principal subspace of a standard
module for an untwisted affine Lie algebra of type $A$ given in
\cite{FS1}--\cite{FS2}:

\begin{definition} \label{def-ps}
\rm Consider any finite-dimensional semisimple Lie algebra and its
associated affine Lie algebra. The {\em principal subspace} of a
highest weight module $V$ for the affine Lie algebra is
$U(\bar{\goth{n}}) \cdot v \subset V$, where $\goth{n}$ is the
nilradical of a a fixed Borel subalgebra of the finite-dimensional Lie
algebra, $\bar{\goth{n}}=\goth{n} \otimes \mathbb{C}[t, t^{-1}]$ and
$v$ is a highest weight vector of $V$.
\end{definition}

In particular, the principal subspace $W(\Lambda_{k,i})$ of
$L(\Lambda_{k,i})$ is
\begin{equation} \label{ps}
W(\Lambda_{k,i})= U(\bar{\goth{n}}) \cdot
v_{\Lambda_{k,i}}
\end{equation}
for $i=0,\dots k$, as in \cite{FS1}.

We have
\begin{equation} \label{zero}
W(\Lambda_{k,i})=U(\bar{\goth{n}}_{-}) \cdot
v_{\Lambda_{k, i}}.
\end{equation}
Set
\begin{equation}
W(\Lambda_{k,k})' = U(\bar{\goth{n}}_{\leq -2}) \cdot v_{\Lambda_{k,k}},
\end{equation}
generalizing (2.8) in \cite{CalLM1}.  Since
$x_{\alpha}(-1)\cdot v_{\Lambda_{k,k}}=0$, we have
\begin{equation} \label{prime}
W(\Lambda_{k,k})' = W(\Lambda_{k,k}),
\end{equation}
generalizing (2.9) in \cite{CalLM1}.

For $i=0,\dots k$, consider the surjective maps
\begin{eqnarray} \label{surj1}
F_{\Lambda_{k,i}}: U(\widehat{\goth{g}}) & \longrightarrow &
L(\Lambda_{k,i}) \\ a &\mapsto& a \cdot v_{\Lambda_{k,i}} .\nonumber
\end{eqnarray}
Restrict $F_{\Lambda_{k,i}}$ to $U(\bar{\goth{n}}_{-})$ and
$F_{\Lambda_{k,k}}$ to $U(\bar{\goth{n}}_{\leq -2})$ and denote these
(surjective) restrictions by $f_{\Lambda_{k,i}}$ and
$f_{\Lambda_{k,k}}'$:
\begin{eqnarray} \label{surj2}
f_{\Lambda_{k,i}}: U(\bar{\goth{n}}_{-}) & \longrightarrow &
W(\Lambda_{k,i})\\ a & \mapsto & a \cdot v_{\Lambda_{k,i}}, \nonumber
\end{eqnarray}
\begin{eqnarray}
f_{\Lambda_{k,k}}': U(\bar{\goth{n}}_{\leq -2}) & \longrightarrow &
W(\Lambda_{k,k})' \\ a & \mapsto & a \cdot v_{\Lambda_{k,
k}}, \nonumber
\end{eqnarray}
generalizing (2.11) and (2.12) in \cite{CalLM1}.  Our main goal is to
give a precise description of the kernels $\mbox{Ker} \;
f_{\Lambda_{k,i}}$ and $\mbox{Ker} \; f_{\Lambda_{k,k}}'$.

For every $t \in \mathbb{Z}$ we consider
the following formal infinite sums:
\begin{equation} \label{R_t}
R_{k,t}=\sum_{m_1+\cdots+m_{k+1}=-t}
x_{\alpha}(m_1)\cdots x_{\alpha}(m_{k+1}).
\end{equation}
For each $t$, $R_{k,t}$ acts naturally on any highest weight
$\widehat{\goth{g}}$-module and, in particular, on each
$L(\Lambda_{k,i})$ for $0 \leq i \leq k$.  For $k=1$ these are the
formal sums $R_t$ introduced in \cite{CalLM1}.

Continuing to generalize the corresponding objects in \cite{CalLM1},
in order to describe $\mbox{Ker} \;
f_{\Lambda_{k,i}}$ and $\mbox{Ker} \; f_{\Lambda_{k,k}}'$ we shall
truncate each $R_{k,t}$ as follows:
\begin{equation} \label{R^0_t}
R_{k,t}^0 =\sum_{{\tiny \begin{array}{c} m_1, \dots, m_{k+1}
\leq -1, \\ m_1+\cdots+m_{k+1}=-t
\end{array}}}x_{\alpha}(m_1)\cdots x_{\alpha}(m_{k+1}), \;
\; t \geq k+1.
\end{equation}
Just as in \cite{CalLM1}, we shall often be viewing $R_{k,t}^0$ as an
element of $U(\bar{\goth{n}})$, and in fact of
$U(\bar{\goth{n}}_{-})$, rather than as an endomorphism of a
$\widehat{\goth{g}}$-module.
In order to describe $\mbox{Ker} \; f'_{\Lambda_{k,k}}$ it will also be
convenient to take $m_1, \dots ,m_{k+1} \leq -2$ in (\ref{R_t}), to
obtain other elements of $U(\bar{\goth{n}})$, which we denote by
$R_{k,t}^1$:
\begin{equation} \label{R^1_t}
R_{k,t}^1=\sum_{{\tiny \begin{array}{c} m_1, \dots ,m_{k+1} \leq
-2,
\\ m_1+\cdots+m_{k+1}=-t \end{array}}}
x_{\alpha}(m_1) \cdots x_\alpha(m_{k+1}), \; \; t
\geq 2(k+1).
\end{equation}

Again as in \cite{CalLM1}, one can view $U(\bar{\goth{n}}_{-})$ and
$U(\bar{\goth{n}}_{\leq -2})$ as the polynomial algebras
\begin{equation}
U(\bar{\goth{n}}_{-})= \mathbb{C}[x_{\alpha}(-1), x_{\alpha}(-2),
\dots ]
\end{equation}
and
\begin{equation}
U(\bar{\goth{n}}_{\leq-2})= \mathbb{C}[x_{\alpha}(-2), x_{\alpha}(-3),
\dots ],
\end{equation}
so that
\begin{equation} \label{deco}
U(\bar{\goth{n}}_{-})= U(\bar{\goth{n}}_{\leq -2}) \oplus
U(\bar{\goth{n}}_{-}) x_{\alpha}(-1)
\end{equation}
and we have the corresponding projection
\begin{equation} \label{rho}
\rho: U(\bar{\goth{n}}_{-}) \longrightarrow U(\bar{\goth{n}}_{\leq
  -2}).
\end{equation}
{}From (\ref{R^0_t}) and (\ref{R^1_t}) we have
\begin{equation}
R_{k,t}^1= \rho (R_{k,t}^0).
\end{equation}
(For $t<2(k+1)$, $R_{k,t}^1=0$.)

Generalizing the corresponding constructions in \cite{CalLM1}, we set
\begin{equation} \label{ideal1}
I_{\Lambda_{k,0}}= \sum_{t \geq k+1} U(\bar{\goth{n}}_{-})R_{k,t}^0
\subset U(\bar{\goth{n}}_{-}),
\end{equation}
\begin{equation} \label{ideal2}
I_{\Lambda_{k,i}}=\sum_{t \geq k+1} U(\bar{\goth{n}}_{-})R_{k,t}^0 +
U(\bar{\goth{n}}_{-})x_{\alpha}(-1)^{k-i+1} \subset
U(\bar{\goth{n}}_{-}) \; \; \mbox{for} \; \; i \geq 0
\end{equation}
(note that (\ref{ideal2}) indeed agrees with (\ref{ideal1})
for $i=0$, since $R_{k,k+1}^0 = x_{\alpha}(-1)^{k+1}$) and
\begin{equation} \label{I'}
I_{\Lambda_{k,k}}'=\sum_{t \geq 2(k+1)} U(\bar{\goth{n}}_{\leq
  -2})R_{k,t}^1 \subset U(\bar{\goth{n}}_{\leq -2}).
\end{equation}

\begin{remark} \label{inclusions}
\em We have the inclusions
\begin{equation} \label{eq0}
I_{\Lambda_{k,0}} \subset I_{\Lambda_{k,1}} \subset \cdots \subset
I_{\Lambda_{k, k-1}} \subset I_{\Lambda_{k,k}}
\end{equation}
among the $U(\bar{\goth{n}}_-)$-ideals $I_{\Lambda_{k,i}}$.  We
also have
\begin{equation} \label{0ik}
I_{\Lambda_{k,i}}=I_{\Lambda_{k,0}} + U(\bar{\goth{n}}_{-})
x_{\alpha}(-1)^{k-i+1} \; \; \mbox{for every} \; \; i \geq 1,
\end{equation}
and this holds for $i=0$ as well.  In addition,
\begin{equation}
\rho  (I_{\Lambda_{k,k}}) = I_{\Lambda_{k,k}}',
\end{equation}
and in fact,
\begin{equation} \label{1-1'}
I_{\Lambda_{k,k}}= I_{\Lambda_{k,k}}' \oplus
U(\bar{\goth{n}}_{-})x_{\alpha}(-1).
\end{equation}
These relations generalize corresponding ones in \cite{CalLM1}.
\end{remark}

\section{Formulations of the main result}
\setcounter{equation}{0}

It is well known that the level $k$ standard
$\widehat{\goth{sl}(2)}$-module $L(\Lambda_{k,0})$ has a natural
vertex operator algebra structure with $v_{\Lambda_{k, 0}}$ as vacuum
vector; the vertex operator map
\begin{eqnarray}
Y(\cdot , x) : L(\Lambda_{k,0}) & \longrightarrow & \mbox{End} \;
L(\Lambda_{k,0}) \; [[x, x^{-1}]] \\ v & \mapsto & Y(v, x)= \sum_{m \in
\mathbb{Z}} v_m x^{-m-1} \nonumber
\end{eqnarray}
has the property
\begin{equation} \label{vertex}
Y(x_{\alpha}(-1)\cdot v_{\Lambda_{k, 0}}, x)=
\sum_{m \in \mathbb{Z}} x_{\alpha}(m) x^{-m-1}.
\end{equation}
It is also well known that each $L(\ki)$, $0 \leq i \leq k$, has a
natural $L(\Lambda_{k,0})$-module structure, with (\ref{vertex})
remaining valid for the module action.  (See \cite{FZ}, \cite{DL},
\cite{Li}, \cite{LL}.)

The standard action of the Virasoro algebra operator $L(0)$ (not to be
confused with the trivial $\widehat{\goth{g}}$-module) provides the
usual grading by {\it conformal weights} on the spaces $L(\ki)$.  We
have
\begin{equation} \label{wt}
\mbox{wt} \; x_{\alpha}(m)= -m
\end{equation}
for $m \in \mathbb{Z}$, where $x_{\alpha}(m)$ is viewed as either an
operator or as an element of $U(\bar{\goth{n}})$.
For any $i$ with $0 \leq i \leq k$,
\begin{equation} \label{form}
\mbox{wt} \; v_{\Lambda_{k,i}} \; = \; \frac{\langle
i\alpha/2,i\alpha/2 + \alpha \rangle}{2(k+2)} \;= \;
\frac{i^2+2i}{4(k+2)}
\end{equation}
(cf. \cite{K}, \cite{DL}, \cite{LL}).

There is also a grading by {\it charge}, given by the eigenvalues of
the operator $\frac{1}{2} \alpha(0)=\frac{1}{2} h(0)$, on
$L(\ki)$. The weight and charge gradings are compatible.  For any $m
\in \mathbb{Z}$, $x_{\alpha}(m)$, viewed as either an operator or as
an element of $U(\bar{\goth{n}})$, has charge $1$. Also,
$v_{\Lambda_{k,i}}$ has charge $\langle i \alpha/2, \alpha/2
\rangle=i/2$.  The principal subspaces $W(\ki)$ are graded by weight
and charge.  For any $m_1, \dots , m_r \in \mathbb{Z}$,
\begin{equation} \label{v0}
x_{\alpha}(m_1) \cdots x_{\alpha}(m_r) \cdot v_{\Lambda_{k,i}} \in
W(\ki)
\end{equation}
has weight $ -m_1-\cdots -m_r+ \frac{i^2+2i}{4(k+2)}$ and charge
$r+\frac{i}{2}$.  See \cite{CLM1}--\cite{CLM2} and \cite{CalLM1} for
further details, background and notation.

\begin{remark} \label{L(0)}
\em
As in \cite{CalLM1}, we have
$$
L(0) \; \mbox{\rm Ker} \; f_{\Lambda_{k,i}} \subset \mbox{\rm Ker}
\; f_{\Lambda_{k,i}} \; \; \mbox{for all} \; \; 0 \leq i \leq k
$$
and
$$
L(0) \; \mbox{Ker} \; f_{\Lambda_{k,k}}' \subset \mbox{Ker} \;
f_{\Lambda_{k,k}}'.
$$
Also, $R_{k,t}^0$ and $R_{k,t}^1$ have conformal weight $t$:
$$
L(0) R_{k,t}^0 = t R_{k,t}^0 \; \; \; \mbox{for all} \; \; \; t
\geq k+1
$$
and
$$ L(0) R_{k,t}^1= t R_{k,t}^1 \; \; \; \mbox{for all} \; \; \; t \geq
2(k+1),
$$
so that in particular, the subspaces $I_{\Lambda_{k,i}}$ and
$I'_{\Lambda_{k,k}}$ are $L(0)$-stable.  Also, $R_{k,t}^0$ and
$R_{k,t}^1$ have charge $k+1$, and the spaces $\mbox{Ker} \;
f_{\Lambda_{k,i}}$, $\mbox{Ker} \; f_{\Lambda_{k,k}}'$,
$I_{\Lambda_{k,i}}$ and $I_{\Lambda_{k,k}}'$ are graded by charge.
Thus these spaces are graded by both weight and charge, and the two
gradings are compatible.
\end{remark}

We will prove the following description of the kernels $\mbox{Ker}\;
f_{\Lambda_{k,i}}$ and $\mbox{Ker} \; f_{\Lambda_{k,k}}'$ (recall
(\ref{deco}) and (\ref{1-1'})):
\begin{theorem} \label{th1}
For any $i = 0, \dots, k$, we have
\begin{equation}
\mbox{ \rm Ker} \; f_{\Lambda_{k,i}}= I_{\Lambda_{k,i}}.
\end{equation}
In particular,
$$
{\rm Ker} \; f'_{\Lambda_{k,k}}= I'_{\Lambda_{k,k}}.
$$
\end{theorem}

As in \cite{CalLM1}, we will actually prove a restatement of this
assertion (see Theorem \ref{th2} below) that uses generalized
Verma modules, in the sense of \cite{L1}, \cite{GL} and \cite{L2},
for $\widehat{\goth{g}}$, and the principal subspaces of these
generalized Verma modules.

The generalized Verma module $N(\Lambda_{k,0})$ is defined as the induced
$\widehat{\goth{g}}$-module
\begin{equation}
N(\Lambda_{k,0}) = U(\widehat{\goth{g}}) \otimes
_{U(\widehat{\goth{g}}_{\geq 0})} \mathbb{C}v_{ \Lambda_{k, 0}}^N,
\end{equation}
where $\goth{g} \otimes \mathbb{C}[t]$ acts trivially and ${\textbf
k}$ acts as the scalar $k$ on $\mathbb{C}v_{\Lambda_{k,0}}^N$;
$v_{\Lambda_{k,0}}^N$ is a highest weight vector.  From the
Poincar\'e-Birkhoff-Witt theorem we have
\begin{equation} \label{pbw}
N(\Lambda_{k,0}) =U(\widehat{\goth{g}}_{<0}) \otimes_{\mathbb{C}}
U(\widehat{\goth{g}}_{\geq 0}) \otimes _{U(\widehat{\goth{g}}_{\geq
0})} \mathbb{C}v_{\Lambda_{k,0}}^N= U(\widehat{\goth{g}}_{<0}) \otimes
_{\mathbb{C}} \mathbb{C}v_{k \Lambda_0}^N =U(\widehat{\goth{g}}_{<0}),
\end{equation}
with the natural identifications.  We similarly define the generalized
Verma module
$$
N(\ki)=U(\widehat{\goth{g}}) \otimes_{U(\widehat{\goth{g}}_{\geq
0})} U_i
$$
for $i = 1,\dots k$, where $U_i$ is an $i+1$-dimensional irreducible
$\goth{g}$-module
and where $\goth{g} \otimes t \mathbb{C}[t]$ acts trivially and
${\textbf k}$ acts by $k$.  By the Poincar\'e-Birkhoff-Witt theorem we
have the identifications
$$
N(\ki)= U(\widehat{\goth{g}}_{<0}) \otimes_{\mathbb{C}} U_i.
$$
For $0 \leq i \leq k$ we have the natural surjective
$\widehat{\goth{g}}$-module maps
\begin{eqnarray} \label{starr}
F_{\Lambda_{k,i}}^N: U(\widehat{\goth{g}}) & \longrightarrow & N(\ki)  \\
a & \mapsto & a \cdot v_{\Lambda_{k,i}}^N, \nonumber
\end{eqnarray}
where $v_{\Lambda_{k,i}}^N$ is a highest weight vector of $U_i$
(cf. (\ref{surj1})).

\begin{remark} \label{star}
\rm The restriction of (\ref{starr}) to
$U(\widehat{\goth{g}}_{<0})$ is a
$U(\widehat{\goth{g}}_{<0})$-module isomorphism for $i=0$ and a
$U(\widehat{\goth{g}}_{<0})$-module injection for $i \geq 1$.
\end{remark}

{}From Definition \ref{def-ps}, the $\bar{\goth{n}}$-submodule
\begin{equation}
W^N(\ki)= U(\bar{\goth{n}}) \cdot v_{\Lambda_{k,i}}^N
\end{equation}
of $N(\ki)$ is the principal subspace of the generalized Verma module
$N(\ki)$, generalizing the corresponding structure in \cite{CalLM1}.
We have
\begin{equation}
W^N(\ki) = U(\bar{\goth{n}}_{-}) \cdot v_{\Lambda_{k,i}}^N.
\end{equation}
We also consider the subspace
\begin{equation}
W^N(\Lambda_{k,k})' = U(\bar{\goth{n}}_{\leq -2}) \cdot
v_{\Lambda_{k,k}}^N
\end{equation}
of $W^N(\Lambda_{k,k})$.

\begin{remark} \label{star-pr}
\rm In view of Remark \ref{star}, the restrictions of
$F_{\Lambda_{k,i}}^N$ to $U(\bar{\goth{n}}_{-})$,
\begin{eqnarray} \label{eq1}
U(\bar{\goth{n}}_{-}) & \longrightarrow & W^N(\Lambda_{k,i}) \\
a & \mapsto & a \cdot v_{\Lambda_{k,i}}^N , \nonumber
\end{eqnarray}
are $\bar{\goth{n}}_{-}$-module isomorphisms and the restriction of
$F_{\Lambda_{k,k}}^N$ to $U(\bar{\goth{n}}_{\leq -2})$,
\begin{eqnarray} \label{eq2}
U(\bar{\goth{n}}_{\leq -2}) & \longrightarrow & W^N(\Lambda_{k,k})' \\
a & \mapsto & a \cdot v_{\Lambda_{k,k}}^N , \nonumber
\end{eqnarray}
is an $\bar{\goth{n}}_{\leq -2}$-module isomorphism.
\end{remark}

In particular, by using (\ref{deco}) we have the natural
identifications
\begin{equation}
W^N(\Lambda_{k,k})' \simeq W^N(\Lambda_{k,k})/ U(\bar{\goth{n}}_{-})
x_{\alpha}(-1) \cdot v_{\Lambda_{k,k}}^N \simeq
U(\bar{\goth{n}}_{\leq -2}).
\end{equation}

Consider the natural surjective $\widehat{\goth{g}}$-module maps
\begin{eqnarray}
\Pi_{\Lambda_{k,i}}: N(\ki) & \longrightarrow & L(\ki) \\
a \cdot v_{\Lambda_{k,i}}^N & \mapsto & a \cdot v_{\Lambda_{k,i}}
\nonumber
\end{eqnarray}
for $a \in U(\widehat{\goth{g}})$ and set
\begin{equation}
N^1(\Lambda_{k,i}) =\mbox{Ker} \; \Pi_{\Lambda_{k,i}}.
\end{equation}
The restrictions of
$\Pi_{\Lambda_{k,i}}$ to $W^N(\ki)$ (respectively,
$W^N(\Lambda_{k,k})'$) are $\bar{\goth{n}}$-module (respectively,
$\bar{\goth{n}}_{\leq -2}$-module) surjections:
\begin{equation} \label{pi}
\pi_{\Lambda_{k,i}}: W^N(\ki) \longrightarrow W(\ki)
\end{equation}
for $0 \leq i \leq k$ and
\begin{equation}
\pi_{\Lambda_{k,k}}': W^N(\Lambda_{k,k})' \longrightarrow
W(\Lambda_{k,k})
\end{equation}
(recall (\ref{prime})).

As in the case of $L(\ki)$, the generalized Verma modules $N(\ki)$ are
compatibly graded by conformal weight and by charge.  We shall
restrict these gradings to the principal subspaces $W^N(\ki)$. The
elements of $W^N(\ki)$ given by (\ref{v0}) with $v_{\Lambda_{k,i}}$
replaced by $v_{\Lambda_{k,i}}^N$ have the same weights and charges as
in those cases.

\begin{remark} \label{L(0)k}
\rm Since the maps $\pi_{\Lambda_{k,i}}$, and $\pi_{\Lambda_{k,k}}'$
commute with the actions of $L(0)$, the kernels $\mbox{Ker} \;
\pi_{\Lambda_{k,i}}$ and $\mbox{Ker} \; \pi_{\Lambda_{k,k}}'$ are
$L(0)$-stable.  These maps also preserve charge, so that $\mbox{Ker}
\; \pi_{\Lambda_{k,i}}$ and $\mbox{Ker} \; \pi_{\Lambda_{k,k}}'$ are
also graded by charge.
\end{remark}

Using Remark \ref{star}, we see that Theorem \ref{th1} can be
reformulated as follows:
\begin{theorem} \label{th2}
For $i = 0,\dots, k$, we have
\begin{equation} \label{eqone}
{\rm Ker} \; \pi_{\Lambda_{k,i}} = I_{\Lambda_{k,i}} \cdot
v_{\Lambda_{k,i}}^N \; \; (\subset N^1(\ki)).
\end{equation}
In particular,
$$
{\rm Ker} \; \pi_{\Lambda_{k,k}}' = I_{\Lambda_{k,k}}' \cdot
v_{\Lambda_{k,k}}^N \; \; (\subset N^1(\Lambda_{k,k})).
$$
\end{theorem}

\section{ Proof of the main result}
\setcounter{equation}{0}

Using the setting of \cite{CLM1}--\cite{CLM2}, with
$P=\frac{1}{2}\mathbb{Z}\alpha$ the weight lattice of $\goth{sl}(2)$,
we have the space
\begin{equation}\label{VP}
V_P=L(\Lambda_0) \oplus L(\Lambda_1)
\end{equation}
and its vertex operator structure.  We shall use the identifications
\begin{equation} \label{vectk=1}
v_{\Lambda_{1,0}}=1 \in L(\Lambda_0) \ \ {\rm and} \ \
v_{\Lambda_{1,1}}=e^{\alpha/2} \cdot v_{\Lambda_{1,0}} \in L(\Lambda_1)
\end{equation}
as in formula (2.5) in \cite{CalLM1} and Section 2 of \cite{CLM1}
(with $v_{\Lambda_0}=v_{\Lambda_{1,0}}$ and
$v_{\Lambda_1}=v_{\Lambda_{1,1}}$, using our current notation for
highest weight vectors).  We consider
\begin{equation}\label{VPk}
V_P^{\otimes \; k} = V_P \otimes \cdots \otimes V_P
\end{equation}
($k$ times).  For any $k$-tuple $(j_1, \dots, j_k)$ with $j_1, \dots,
j_k \in \{0,1\}$ we consider the element
\begin{equation}\label{vj1jk}
v_{j_1, \dots, j_k}=v_{\Lambda_{1, j_1}} \otimes \cdots \otimes
v_{\Lambda_{1, j_k}} \in V_P^{\otimes k},
\end{equation}
where exactly $k-i$ indices $j_l$ ($l=1, \dots, k$) are equal to 0 (and
exactly $i$ indices are equal to 1); recall (\ref{vectk=1}). This
vector is of course a highest weight vector for
$\widehat{\goth{sl}(2)}$, and
\begin{equation}
L(\Lambda_{k,i}) \simeq U(\widehat{\goth{g}}) \cdot v_{j_1,
\dots, j_k} \subset V_P^{\otimes \; k}
\end{equation}
(cf. \cite{K} and \cite{CLM2}), using the natural extension to
$U(\widehat{\goth{g}})$ of the usual comultiplication
\begin{equation} \label{action}
a \cdot v= (a \otimes \cdots \otimes 1 + \cdots
+ 1 \otimes \cdots \otimes a) v
\end{equation}
for $a \in \widehat{\goth{g}}$ and $v \in V_P^{\otimes \; k}$.  As in
\cite{CLM2} we will use the embeddings
\begin{equation} \label{L}
\iota_{j_1, \dots, j_k} : L(\Lambda_{k,i}) \hookrightarrow
V_P^{\otimes \; k} \; \; \; \mbox{for} \; \; 0 \leq i \leq k,
\end{equation}
uniquely determined by the identifications
\begin{equation} \label{v}
v_{\Lambda_{k,i}}=v_{j_1, \dots, j_k}.
\end{equation}
Of course, this element $v_{\Lambda_{k,i}}$ and the embedding
(\ref{L}) depend on $j_1, \dots, j_k$.

Recall the linear isomorphism (3.20) in \cite{CLM1},
\[
e^{\alpha/2}:V_P \longrightarrow V_P,
\]
and consider the linear isomorphism
\begin{equation} \label{map}
\begin{array}{c} e^{\alpha/2}_{(k)}=\underbrace{e^{\alpha/2}
\otimes \cdots \otimes e^{\alpha/2}} : V_P^{\otimes \; k}
\longrightarrow V_P^{\otimes \; k}. \\ k \; \mbox{times} \ \ \ \ \ \ \ \
\end{array}
\end{equation}
We denote by
$e^{\alpha/2}_{(k,i)}$ the restriction of $e^{\alpha/2}_{(k)}$ to the
principal subspace $W(\ki)$ of $L(\ki)$, using an embedding of the
form (\ref{L}).  The action of $e^{\alpha/2}_{(k)}$ on
$L(\Lambda_{k,i})$ and its action $e^{\alpha/2}_{(k,i)}$ on $W(\ki)$
as well as other features of this map are given as follows (see Lemma
3.2 of \cite{CLM2} and its proof):

\begin{lemma} \cite{CLM2} Fix $j_1, \dots, j_k$ as in (\ref{vj1jk})--(\ref{v})
and consider the standard module $L(\Lambda_{k,i})$ embedded in
$V_P^{\otimes \; k}$ via $\iota_{j_1, \dots, j_k}$. The image of the
restriction of $e^{\alpha/2}_{(k)}$ to $L(\Lambda_{k,i})$ lies in
$L(\Lambda_{k, k-i})$, embedded in $V_P^{\otimes \; k}$ via
$\iota_{1-j_1, \dots, 1-j_k}$.  For any $i$ with $0 \leq i \leq k$, we
have
\begin{equation} \label{shift}
e^{\alpha/2}_{(k,i)} : W(\ki) \longrightarrow
W(\Lambda_{k, k-i}).
\end{equation}
When $i=0$ the map (\ref{shift}) is a linear isomorphism.
We also have
\begin{equation} \label{def-shift}
e^{\alpha/2}_{(k, i)} \; (x_{\alpha}(m_1)\cdots x_{\alpha}(m_r) \cdot
v_{\Lambda_{k,i}})= \frac{1}{i!} x_{\alpha}(m_1-1) \cdots
x_{\alpha}(m_r-1)x_{\alpha}(-1)^i \cdot v_{\Lambda_{k,k-i}}
\end{equation}
for any $m_1, \dots, m_r \in \mathbb{Z}$. $\; \; \; \; \Box$
\end{lemma}

We emphasize that according to our notation, the embeddings of the two
spaces $W(\Lambda_{k,i})$ and $W(\Lambda_{k, k-i})$ in (\ref{shift})
are ``opposite'' even when $i$ and $k-i$ happen to coincide.

We now generalize the lifting procedures in \cite{CalLM1}.
For each $i=0, \dots, k$ we construct a lifting
\begin{equation} \label{hatshift}
\widehat{e^{\alpha/2}_{(k,i)}}: W^N(\ki) \longrightarrow W^N(\Lambda_{k,k-i})
\end{equation}
 of
\begin{equation} \label{hat1}
e^{\alpha/2}_{(k,i)}: W(\ki) \longrightarrow W(\Lambda_{k,k-i}),
\end{equation}
making the diagram
$$ \CD W^N(\ki) @> \widehat{e^{\alpha/2}_{(k,i)}}> >
W^N(\Lambda_{k,k-i}) \\ @V\pi_{\Lambda_{k,i}}VV
@V\pi_{\Lambda_{k,k-i}}VV \\ W(\ki)\ @> {e^{\alpha/2}_{(k,i)}}> >
W(\Lambda_{k,k-i})  \endCD
$$
commute; here, in (\ref{hat1}) we continue to use the particular embeddings
depending on $j_1, \dots, j_k$ used in (\ref{shift}).
In fact, for any $i$ with $0 \leq i \leq k$ and any integers $m_1, \dots,
m_r <0$ we set
\begin{equation} \label{def-hatshift}
\widehat{e^{\alpha/2}_{(k, i)}} \; ( x_{\alpha}(m_1) \cdots
x_{\alpha}(m_r) \cdot v_{\Lambda_{k,i}}^N)= \frac{1}{i!}x_{\alpha}(m_1-1)
\cdots x_{\alpha}(m_r-1)x_{\alpha}(-1)^i \cdot
v_{\Lambda_{k,k-i}}^N,
\end{equation}
which is well defined, since $U(\bar{\goth{n}}_{-})$, viewed as
the polynomial algebra
$$\mathbb{C}[x_{\alpha}(-1), x_{\alpha}(-2), \dots],$$ maps
isomorphically onto $W^N(\ki)$ under the map (\ref{eq1}). This gives
our desired lifting (\ref{hatshift}).

For the case $i=0$, the map (\ref{hatshift}) is a linear isomorphism
onto the subspace $W^N(\Lambda_{k,k})'$:
\begin{equation}\label{shiftiso}
\widehat{e^{\alpha/2}_{(k, 0)}}: W^N(\Lambda_{k,0}) \longrightarrow
W^N(\Lambda_{k,k})',
\end{equation}
a lifting of the linear isomorphism
\begin{equation}\label{hat2}
e^{\alpha/2}_{(k,0)}: W(\Lambda_{k,0}) \longrightarrow W(\Lambda_{k,k});
\end{equation}
the diagram
$$ \CD W^N(\Lambda_{k,0}) @> \widehat{e^{\alpha/2}_{(k,0)}}> \sim >
W^N(\Lambda_{k,k})' \\ \ @V\pi_{\Lambda_{k,0}}VV @V\pi'_{\Lambda_{k,k}}VV
\\ W(\Lambda_{k,0})\ @> {e^{\alpha/2}_{(k,0)}}>\sim > W(\Lambda_{k,k})
\endCD
$$
commutes.  Indeed, since
$$
W^N(\Lambda_{k,k})'=U(\bar{\goth{n}}_{\leq -2}) \cdot
v_{\Lambda_{k,k}}^N= \widehat{e^{\alpha/2}_{(k,0)}} \;
(U(\bar{\goth{n}}_{-}) \cdot v_{\Lambda_{k,0}}^N),
$$
the linear map (\ref{shiftiso}) is surjective, and by Remark
\ref{star-pr} it is also injective and thus a linear isomorphism.
Denote by
\begin{equation} \label{shift1}
(\widehat{e^{\alpha/2}_{(k,0)}})^{-1}=\widehat{e^{-\alpha/2}_{(k,0)}}:
W^N(\Lambda_{k,k})' \longrightarrow W^N(\Lambda_{k,0})
\end{equation}
its inverse; the map $\widehat{e^{-\alpha/2}_{(k,0)}}$ is
correspondingly a lifting of the inverse
\begin{equation} \label{shift11}
e^{-\alpha/2}_{(k,0)}: W(\Lambda_{k,k}) \longrightarrow
W(\Lambda_{k,0}).
\end{equation}

\begin{remark}
\rm We have just noticed that, as in the $k=1$ special case of
\cite{CalLM1}, the image of $W^N(\Lambda_{k,0})$ under the map
$\widehat{e^{\alpha/2}_{(k,0)}}$ is the subspace
$W^N(\Lambda_{k,k})' \subset W^N(\Lambda_{k,k})$ and not the full
space $W^N(\Lambda_{k,k})$. Both the maps (\ref{shiftiso}) and
(\ref{hat2}) are isomorphisms, while the map (\ref{hatshift}) for
$i=0$, from $W^N(\Lambda_{k,0})$ to $W^N(\Lambda_{k,k})$, is only
an injection.
\end{remark}

\begin{remark}
\rm The restriction
\begin{equation} \label{k}
\widehat{e^{\alpha/2}_{(k,k)}}:W^N(\Lambda_{k,k})' \longrightarrow
W^N(\Lambda_{k,0})
\end{equation}
of (\ref{hatshift}) for $i=k$ to $W^N(\Lambda_{k,k})'$ is a lifting of
\begin{equation}
e^{\alpha/2}_{(k,k)}:W(\Lambda_{k,k}) \longrightarrow W(\Lambda_{k,0}),
\end{equation}
making the diagram
$$ \CD W^N(\Lambda_{k,k})' @> \widehat{e^{\alpha/2}_{(k,k)}}> >
W^N(\Lambda_{k,0}) \\ \ @V\pi'_{\Lambda_{k,k}}VV @V\pi_{\Lambda_{k,0}}VV
\\ W(\Lambda_{k,k})\ @> {e^{\alpha/2}_{(k,k)}}> > W(\Lambda_{k,0})
\endCD
$$
commute; it is an injection and not a surjection.  The maps
(\ref{shiftiso}) and (\ref{k}) were used in \cite{CalLM1} for $k=1$.
\end{remark}

Now we describe the actions of our liftings (\ref{hatshift}) on the
spaces $I_{\Lambda_{k,i}} \cdot v_{\Lambda_{k,i}}^N$:

\begin{lemma}\label{lemma2} For any $i$ with $0 \leq i \leq k$, we have
\begin{equation} \label{l2}
\widehat{e^{\alpha/2}_{(k,i)}} \; (I_{\Lambda_{k,i}} \cdot
v_{\Lambda_{k,i}}^N) \subset I_{\Lambda_{k,k-i}} \cdot
v_{\Lambda_{k,k-i}}^N.
\end{equation}
\end{lemma}
\noindent {\it Proof:} By (\ref{ideal2}) we have
\begin{equation} \label{desc1}
I_{\Lambda_{k,i}} \cdot v_{\Lambda_{k,i}}^N= \sum_{t \geq k+1}
U(\bar{\goth{n}}_{-}) R_{k, t}^0 \cdot v_{\Lambda_{k,i}}^N +
U(\bar{\goth{n}}_{-}) x_{\alpha}(-1)^{k-i+1} \cdot v_{\Lambda_{k,
i}}^N
\end{equation}
and
\begin{equation}\label{desc2}
I_{\Lambda_{k,k-i}} \cdot v_{\Lambda_{k,k-i}}^N= \sum_{t \geq k+1}
U(\bar{\goth{n}}_{-}) R_{k, t}^0 \cdot v_{\Lambda_{k,k-i}}^N +
U(\bar{\goth{n}}_{-}) x_{\alpha}(-1)^{i+1} \cdot v_{\Lambda_{k,
k-i}}^N.
\end{equation}
We use (\ref{def-hatshift}).  For any $t \geq k+1$,
\begin{eqnarray}
&&\widehat{e^{\alpha/2}_{(k,i)}} \; (R_{k,t}^0 \cdot
v_{\Lambda_{k,i}}^N) \nonumber \\
&=& \frac{1}{i!} \sum_{{ \tiny{\begin{array}{c} m_1, \dots ,m_{k+1} \leq
-1 \\ m_1+\cdots+m_{k+1}=-t \end{array}}}} x_{\alpha}(m_1-1) \cdots
x_\alpha(m_{k+1}-1) x_{\alpha}(-1)^{i} \cdot v_{\Lambda_{k, k-i}}^N
\nonumber \\
&=& \frac{1}{i!} R^1_{k, t+k+1} x_{\alpha}(-1)^i \cdot v_{\Lambda_{k,
k-i}}^N \nonumber \\
&=& \frac{1}{i!} R^0_{k, t+k+1} x_{\alpha}(-1)^i \cdot
v_{\Lambda_{k, k-i}}^N +a x_{\alpha}(-1)^{i+1} \cdot v_{\Lambda_{k,
k-i}}^N \in I_{\Lambda_{k,k-i}} \cdot v_{\Lambda_{k,k-i}}^N, \nonumber
\end{eqnarray}
where $a \in U(\bar{\goth{n}}_-)$. We also have
\begin{eqnarray} \label{trick}
&& \widehat{e^{\alpha/2}_{(k,i)}} (x_{\alpha}(-1)^{k-i+1} \cdot
v_{\Lambda_{k,i}}^N )= \frac{1}{i!} x_{\alpha}(-2)^{k-i+1}
x_{\alpha}(-1)^{i} \cdot v_{\Lambda_{k, k-i}}^N  \\ &&=
\gamma R_{k, 2k-i+2}^0 \cdot v_{\Lambda_{k, k-i}}^N+b
x_{\alpha}(-1)^{i+1}\cdot v_{\Lambda_{k, k-i}}^N \in
I_{\Lambda_{k,k-i}} \cdot v_{\Lambda_{k,k-i}}^N, \nonumber
\end{eqnarray}
where $\gamma$ is a nonzero scalar and $b \in U(\bar{\goth{n}}_-)$.  Indeed,
the expression $R_{k, 2k-i+2}^0$ does not have any terms
involving $x_{\alpha}(-1)^t$ with $0 \leq t <i$, since if there is
such a term $cx_{\alpha}(-1)^t$, with $c \in U(\bar{\goth{n}}_{-})$ a
product of $k+1-t$ elements $x_{\alpha}(m)$, where each $m \leq -2$,
then
$$
\mbox{wt} \; (c x_{\alpha}(-1)^t) \geq 2(k+1-t)+t=2k+2-t > 2k+2-i =
\mbox{wt} \; (R_{k, 2k-i+2}^0),
$$
and this contradicts the fact that $cx_{\alpha}(-1)^t$ is a summand of
$R_{k, 2k-i+2}^0$. We also observe that
$x_{\alpha}(-2)^{k-i+1}x_{\alpha}(-1)^i$ is the only type of term in the sum
$R^0_{2k-i+2}$ involving $x_{\alpha}(-1)^i$. This proves
(\ref{trick}), and hence (\ref{l2}).  $\; \; \; \Box$
\vspace{1em}

\begin{remark}
\rm We have
\begin{equation} \label{0=1}
\widehat{e^{\alpha/2}_{(k,0)}} \; (I_{\Lambda_{k,0}} \cdot
v_{\Lambda_{k,0}}^N)= I_{\Lambda_{k,k}}' \cdot v_{\Lambda_{k,k}}^N.
\end{equation}
Indeed, for any $t \geq k+1$,
$$
\widehat{e^{\alpha/2}_{(k,0)}} \; (R_{k,t}^0 \cdot v_{\Lambda_{k,0}}^N)
= R_{k, t+k+1}^1 \cdot v_{\Lambda_{k,k}}^N,
$$
and from the descriptions (\ref{ideal1}) and (\ref{I'}) of the ideals
$I_{\Lambda_{k,0}}$ and $I_{\Lambda_{k,k}}'$ we see that (\ref{0=1})
holds.  The $k=1$ case of (\ref{0=1}) was used in \cite{CalLM1}.
\end{remark}

For the reader's convenience we recall from \cite{CalLM1} the shift, or
translation, automorphism
\begin{equation} \label{translation} \tau:
U(\bar{\goth{n}}) \longrightarrow U(\bar{\goth{n}})
\end{equation}
given by
$$
\tau (x_{\alpha}(m_1) \cdots x_{\alpha}(m_r))= x_{\alpha}(m_1-1)
\cdots x_{\alpha}(m_r-1)
$$
for any integers $m_1, \dots, m_k$.
For any integer $s$, the
$s^{\mbox{th}}$ power
\begin{equation}
 \tau ^s: U(\bar{\goth{n}}) \longrightarrow U(\bar{\goth{n}})
\end{equation}
is given by
$$
\tau ^s (x_{\alpha}(m_1) \cdots x_{\alpha}(m_r))=
x_{\alpha}(m_1-s) \cdots x_{\alpha}(m_r-s).
$$
Recall from Remark 3.1 in \cite{CalLM1} that for any nonzero element
$a \in U(\bar{\goth{n}})$ homogeneous with respect to both the weight
and charge gradings such that $a$ has positive charge, the element
$\tau^s(a)$ has the same properties, and
\begin{equation} \label{greater}
\mbox{wt} \; \tau^s (a) > \mbox{wt} \; a \; \; \; \mbox{for} \; \;
\; s > 0
\end{equation}
and
\begin{equation} \label{smaller}
\mbox{wt} \; \tau^s (a) < \mbox{wt} \; a \; \; \; \mbox{for} \; \;
\; s <0.
\end{equation}
If $a$ is a constant, that is, $a$ has charge zero,
then
\begin{equation} \label{constant}
\tau^s(a)= a
\end{equation}
and $\tau^s(a)$ and $a$ have the same weight  and charge.

Generalizing Remark 3.2 in \cite{CalLM1}, we now have:

\begin{remark} \label{rem}
\rm Using the map $\tau$ we can re-express (\ref{shift}) and
(\ref{hatshift}) as follows:
\begin{equation} \label{t1}
e^{\alpha/2}_{(k,i)} (a \cdot v_{\Lambda_{k,i}})=\frac{1}{i!}
\tau(a)x_\alpha(-1)^i \cdot v_{\Lambda_{k,k-i}}, \; \; a \in
U(\bar{\goth{n}}), \; \; 0 \leq i \leq k
\end{equation}
and
\begin{equation} \label{t2}
\widehat{e^{\alpha/2}_{(k,i)}} ( a \cdot v_{\Lambda_{k,i}}^N) =
\frac{1}{i!}  \tau (a) x_\alpha(-1)^i \cdot v_{\Lambda_{k,k-i}}^N, \;
\; a \in U(\bar{\goth{n}}_{-}), \; \; 0 \leq i \leq k.
\end{equation}
(recall (\ref{def-shift}) and (\ref{def-hatshift})).
\end{remark}

Lemma 3.3 in \cite{CalLM1} generalizes to:

\begin{lemma} \label{small}
We have
$$
\tau ( I_{\Lambda_{k,0}}) \subset I_{\Lambda_{k,0}} +
U(\bar{\goth{n}}_{-}) x_{\alpha}(-1)=I_{\Lambda_{k,k}}.
$$
\end{lemma}
\noindent {\it Proof:} Let $t \geq k+1$. Then
$$
\tau (R_{k,t}^0)= R_{k,t+k+1}^1 = R_{k,t+k+1}^0 + a \cdot x_{\alpha}(-1)
$$
for some $a \in U(\bar{\goth{n}}_-)$. $\; \; \; \; \Box$
\vspace{1em}

Intertwining vertex operators (in the sense of \cite{FHL} and
\cite{DL}) among triples of $L(\Lambda_{k,0})$-modules play an
important role in this paper, as they did in \cite{CLM1}--\cite{CLM2}
and \cite{CalLM1}. The following theorem is well known:

\begin{theorem} \cite{FZ} \label{fusions}
For integers $i,j$ and $m$ with $0 \leq i,j, m
\leq k$, write
$$I {L(\Lambda_{k,m}) \choose L(\ki) \ \ L(\Lambda_{k,j})}$$
for the vector space of intertwining operators of type
${L(\Lambda_{k,m}) \choose L(\ki) \ \ L(\Lambda_{k,j})}.$
The dimensions of these spaces (the fusion rules) are given by:
$${\rm dim} \ I {L(\Lambda_{k,m})
\choose L(\ki) \ \ L(\Lambda_{k,j})}=1$$ if and only if
$$i+j-2 \, {\rm max} \{0,i+j-k\} \geq m \geq i+j - 2 \, {\rm min} \{i,j\}, \
m \equiv i+j \ {\rm mod} \ 2,$$
and otherwise, the fusion rule is zero.  $\; \; \; \Box$
\end{theorem}
As a consequence, we see that $L(\Lambda_{k,0})$ and
$L(\Lambda_{k,k})$ are ``group-like'' elements in the fusion ring at
level $k$. Such modules are sometimes called ``simple currents.''

An intertwining operator $\mathcal{Y}( \cdot, x)$ of
type ${L(\Lambda_{k,m}) \choose L(\ki) \ \ L(\Lambda_{k,j})}$
satisfies the condition
$$
\mathcal{Y} (v_{\Lambda_{k,i}}, x) \in x^{{\rm{wt}}\;
v_{\Lambda_{k,m}}-{\rm{wt}} \; v_{\Lambda_{k, i}}-{\rm{wt}} \;
v_{\Lambda_{k, j}}} \mbox{Hom} \; (L(\Lambda_{k,j}), L(\Lambda_{k,m}))
[[ x, x^{-1}]]
$$
(cf. \cite{FHL}, \cite{FZ}). Denote by
$$
{\mathcal Y}_c(v_{\Lambda_{k,i}}, x) \in \mbox{Hom}\;
(L(\Lambda_{k,j}), L(\Lambda_{k,m}))
$$
the constant term of $x^{-{\rm{wt}}\; v_{\Lambda_{k,m}}+{\rm{wt}} \;
v_{\Lambda_{k, i}}+{\rm{wt}} \; v_{\Lambda_{k,j}}}{\mathcal
Y}(v_{\Lambda_{k,i}}, x)$. Then as in
\cite{CLM2}, we have
\begin{equation} \label{p1}
[x_{\alpha}(n), {\mathcal Y}(v_{\Lambda_{k,i}}, x)]=0 \; \;
\mbox{for all} \; \; n \in \mathbb{Z}.
\end{equation}
If ${\mathcal Y}_c (v_{\Lambda_{k,i}} , x)v_{\Lambda_{k,j}}$ is
nonzero then it is a highest weight vector of $L(\Lambda_{k,m})$, so
that
\begin{equation} \label{p2}
{\mathcal Y}_c (v_{\Lambda_{k,i}} , x)v_{\Lambda_{k,j}}=
\gamma v_{\Lambda_{k,m}},
\end{equation}
where $\gamma \neq 0$; this will hold for our cases below.
Using these remarks about intertwining operators and constant
terms we prove the following:

\begin{lemma} \label{kerns}
For any $i$ with $0 \leq i <k$ we have
\begin{equation} \label{incl1}
\mbox{\rm Ker} \; f_{\Lambda_{k,i}} \subset \mbox{\rm Ker} \; f_{\Lambda_{k,
i+1}},
\end{equation}
so that
\begin{equation} \label{incl2}
\mbox{\rm Ker} \; f_{\Lambda_{k,0}} \subset \mbox{\rm Ker} \;
f_{\Lambda_{k,1}} \subset \cdots \subset \mbox{\rm Ker} \;
f_{\Lambda_{k,k}}.
\end{equation}
\end{lemma}
\noindent
{\it Proof:} We consider a nonzero
intertwining operator ${\mathcal Y}$ of type
$$
{L(\Lambda_{k,i+1}) \choose L(\Lambda_{k,1}) \ \
L(\Lambda_{k,i})};
$$
the corresponding fusion rule is one by Theorem
\ref{fusions}. Consider ${\mathcal Y}_c(v_{\Lambda_{k,1}}, x)$, the
constant term of the nonzero operator $x^{-{\rm{wt}} \;
v_{\Lambda_{k,i+1}}+{\rm{wt}} \; v_{\Lambda_{k,1}} + {\rm{wt}} \;
v_{\Lambda_{k,i}}}{\mathcal Y}( v_{\Lambda_{k,1}}, x)$.

Let $a \in U(\bar{\goth{n}}_{-})$ be such that $a \in \mbox{Ker} \
f_{\Lambda_{k,i}}$, so that $a \cdot v_{\Lambda_{k,i}}=0$. By applying
the map ${\mathcal Y}_c(v_{\Lambda_{k,1}}, x)$ to $a \cdot
v_{\Lambda_{k,i}}$ and using (\ref{p1}) and (\ref{p2})
we obtain
$$
\gamma a \cdot v_{\Lambda_{k, i+1}}=0 \; \; \mbox{with} \; \; \gamma
\neq 0,
$$
so that $$a \in \mbox{Ker} \; f_{\Lambda_{k,i+1}},
$$
as desired.  $\: \: \Box$
\vspace{1em}

\begin{remark}
\rm The maps ${\mathcal Y}_c(v_{\Lambda_{k,1}}, x)$ $(0 \leq i <k)$
used here are exactly the same as the constant-term maps crucially
used in Theorem 4.2 (formula (4.44)) of \cite{CLM2}.
\end{remark}

\begin{remark} \rm In order to construct $\mathcal{Y}_c$ and prove Lemma
\ref{kerns} we do not in fact need results from \cite{FZ}.  The
construction of $\mathcal{Y}_c$ follows easily from results in Chapter
13 of \cite{DL}, while Lemma \ref{kerns} follows from the relation
${\rm Ker} \ f_{\Lambda_{1,0}} \subset {\rm Ker} \ f_{\Lambda_{1,1}}$
and (\ref{L}) (cf. also Chapter 13 of \cite{DL}).
\end{remark}

Our next goal is to prove the main result, Theorem \ref{th1}, or
equivalently, Theorem \ref{th2} (formula (\ref{eqone})), which is what
we will in fact prove.

We notice first the inclusion
\begin{equation} \label{first inclusions}
I_{\Lambda_{k,i}} \cdot v_{\Lambda_{k,i}}^N \subset \mbox{Ker} \;
\pi_{\Lambda_{k,i}}, \; \; 0 \leq i \leq k.
\end{equation}
Indeed, as is well known, the $(k+1)$-st power of the vertex operator
$Y(x_{\alpha}(-1) \cdot v_{\Lambda_{k,0}}, x)$ is well defined (the
components $x_{\alpha}(m)$, $m \in \mathbb{Z}$, of this vertex
operator commute) and equals zero on each $L(\ki)$, and in particular
on $W(\ki)$.  The expansion coefficients of
$Y(x_{\alpha}(-1) \cdot v_{\Lambda_{k, 0}}, x)^{k+1}$ are the
operators $R_{k,-t}$, $t \in \mathbb{Z}$:
\begin{equation}
Y(x_{\alpha}(-1) \cdot v_{\Lambda_{k,0}}, x)^{k+1}=
\sum_{t \in \mathbb{Z}} \left (
\sum_{m_1+m_2+\cdots +m_{k+1}=t} x_{\alpha}(m_1)x_{\alpha}(m_2) \cdots
x_\alpha(m_{k+1}) \right )x^{-t-k-1}
\end{equation}
(recall (\ref{R_t}) and (\ref{vertex})).  Thus the operators
(\ref{R^0_t}) annihilate the highest weight vector
$v_{\Lambda_{k,i}},$ and (\ref{first inclusions}) follows.

Before we prove our main result for the general level $k \geq 1$
(Theorem \ref{th2}) we first prove this result for $k=1$, for the
reasons mentioned in Remark \ref{reasons} below.
We have $i=0,1$, and we shall use the notation
$\Lambda_0$ and $\Lambda_1$ instead of $\Lambda_{1,0}$ and
$\Lambda_{1,1}$ (recall (\ref{Lambdakinotation})).

{\it Proof of the $k=1$ case of Theorem \ref{th2}:}
By (\ref{first inclusions}) it is sufficient to show that
\begin{equation}\label{sufftoshow}
\mbox{Ker} \; \pi_{\Lambda_{i}} \subset I_{\Lambda_{i}} \cdot
v_{\Lambda_{i}}^N \; \; \mbox{for} \; \; i=0,1.
\end{equation}
We will prove this by contradiction. Assume then that there exists $a
\in U(\bar{\goth{n}}_{-})$ such that
\begin{equation} \label{form1}
a \cdot v_{\Lambda_i}^N \in \mbox{Ker} \; \pi_{\Lambda_i} \; \;
\mbox{but} \; \; a \cdot v_{\Lambda_i}^N \notin I_{\Lambda_i} \cdot
v_{\Lambda_i}^N \; \; \mbox{for}\; \; i=0 \; \; \mbox{or} \; \; 1.
\end{equation}
By Remarks \ref{L(0)} and \ref{L(0)k} we may and do assume that
$a$ is doubly homogeneous, that is, homogeneous with respect to
the weight and charge gradings.  By the second statement in
(\ref{form1}), $a$ is nonzero, and by the first statement in
(\ref{form1}), $a$ is in fact nonconstant, so that $a$ has
positive weight and positive charge.  Let
\begin{equation}\label{def-L}
L={\rm min} \{ {\rm wt} \; d \; | \; d \in U(\bar{\goth{n}}_-) \;
{\rm{doubly}}\;{\rm{ homogeneous}}\; {\rm{such \ that}} \ (\ref{form1}) \
{\rm{holds \ for}}\; d \}.
\end{equation}
Any such element $d$ is nonzero and in fact nonconstant (just as for
the chosen element $a$), so that any such $d$ has positive weight and
charge; thus $L > 0$.  We further assume that $\mbox{wt} \; a
=L$. Note that $i$ might be $0$ or $1$ or both. We shall show that in
fact $i$ cannot be $1$, and then we shall use this to show that $i$
cannot be $0$, giving our desired contradiction.

By (\ref{deco}) we have a unique decomposition
\begin{equation} \label{form2}
a=r_0x_{\alpha}(-1)+s_0
\end{equation}
with $r_0 \in U(\bar{\goth{n}}_{-})$ and $s_0 \in
U(\bar{\goth{n}}_{\leq -2})$. The elements $r_0$ and $s_0$ are doubly
homogeneous, and in
fact,
\begin{equation} \label{weights-1}
\mbox{wt} \; r_0 = \mbox{wt} \; a -1, \; \; \; \mbox{wt} \; s_0 =
\mbox{wt} \; a;
\end{equation}
similarly, the charge of $r_0$ is one less than that of $a$ and the
charges of $s_0$ and $a$ are equal.  Applying $\tau^{-1}$ to
(\ref{form2}) gives
\begin{equation} \label{form22}
\tau^{-1} (a)= \tau^{-1} (r_0) x_{\alpha}(0)+ \tau^{-1} (s_0),
\end{equation}
and $\tau^{-1} (s_0)$ is doubly homogeneous,
\begin{equation} \label{form222}
\tau^{-1} (s_0) \in U(\bar{\goth{n}}_{-}),
\end{equation}
and
\begin{equation} \label{wttau-1v0}
\mbox{wt} \; \tau^{-1} (s_0) < \mbox{wt} \; a,
\end{equation}
{}from (\ref{smaller}) and the fact that the charge of $a$ and hence of
$s_0$ is positive.

Suppose now that $i=1$. Then we have
\begin{equation} \label{form3}
a \cdot v_{\Lambda_{1}}^N \in \mbox{Ker} \; \pi_{\Lambda_{1}}
\; \; \; \mbox{but} \; \; \; a \cdot v_{\Lambda_{1}}^N \notin
I_{\Lambda_{1}} \cdot v_{\Lambda_{1}}^N \; \; \; (\mbox{that is,} \; \;
\; a \notin I_{\Lambda_{1}}),
\end{equation}
where $a$ is doubly homogeneous and $\mbox{wt} \; a =L$ (recall
(\ref{def-L})). We are going to show that there exists a doubly
homogeneous element of $U(\bar{\goth{n}}_{-})$, namely, $\tau^{-1}
(s_0)$, whose weight is less than $L$ and which satisfies
(\ref{form1}).  We note that $s_0 \neq 0$, because $a \notin
U(\bar{\goth{n}}_{-})x_{\alpha}(-1)$, by (\ref{0ik}) and
(\ref{form3}).  We have seen that $\tau^{-1} (s_0)$ is doubly
homogeneous and that its weight is less than $\mbox{wt} \; a$.
Since $s_0 \cdot v_{\Lambda_{1}}^N \in \ {\rm Ker} \
\pi_{\Lambda_{1}}$ (by (\ref{form2}) and (\ref{form3})), we have $s_0
\cdot v_{\Lambda_{1}}=0$. We also have
$$e^{\alpha/2}_{(1,0)}(\tau^{-1}(s_0) \cdot v_{\Lambda_{0}})= s_0
\cdot v_{\Lambda_{1}}= 0$$
(recall (\ref{t1})), which together with the injectivity of
$e^{\alpha/2}_{(1,0)}$ implies
\be
\label{form4}
\tau^{-1}(s_0) \cdot v_{\Lambda_{0}}^N \in {\rm
Ker} \; \pi_{\Lambda_{0}}.
\ee
We also have
\begin{equation} \label{form5}
\tau^{-1} (s_0) \cdot v_{\Lambda_{0}}^N \notin I_{\Lambda_{0}}
\cdot v_{\Lambda_{0}}^N.
\end{equation}
Indeed, if (\ref{form5}) does not hold, then $\tau^{-1} (s_0) \in
I_{\Lambda_{0}}$, and by Lemma \ref{small} we get $s_0 \in
I_{\Lambda_{1}}$. Now (\ref{form2}) yields $ a \in I_{\Lambda_{1}}$,
and thus $a \cdot v_{\Lambda_1}^N \in I_{\Lambda_1} \cdot
v_{\Lambda_1}^N$, contradicting (\ref{form3}). Hence (\ref{form5})
holds.  Now (\ref{form4}) and (\ref{form5}) give a contradiction since
$\tau^{-1} (s_0)$ is a doubly homogeneous element satisfying
(\ref{form1}) but whose weight is less than $\mbox{wt} \; a = L$.  We
have shown that $i$ cannot be $1$.

Now we may and do assume that $i=0$, that is,
\begin{equation} \label{form55}
a \cdot v_{\Lambda_0}^N \in \mbox{Ker} \; \pi_{\Lambda_0} \; \;
\mbox{but} \; \; a \cdot v_{\Lambda_0}^N \notin I_{\Lambda_0} \cdot
v_{\Lambda_0}^N,
\end{equation}
where $a$ is doubly homogeneous of weight $L$ (recall
(\ref{def-L})). Since $a \cdot v_{\Lambda_0}^N \in \mbox{Ker} \;
\pi_{\Lambda_0}$,
$$
a \cdot v_{\Lambda_0} =0 \; \; \mbox{in} \; \; W(\Lambda_0)
$$
and by Lemma \ref{kerns} we obtain
\begin{equation} \label{form6}
a \cdot v_{\Lambda_1}=0 \; \; \mbox{in} \; \; W(\Lambda_1).
\end{equation}
Hence $a \cdot v_{\Lambda_1}^N \in \mbox{Ker} \; \pi_{\Lambda_1}$, and
so by what we have just proved (that $i$ cannot be $1$), we
obtain
$$
a \cdot v_{\Lambda_1}^N \in I_{\Lambda_1} \cdot v_{\Lambda_1}^N,
$$
and so
$$a \in I_{\Lambda_1}.$$
Our goal is to show that in fact $a \in I_{\Lambda_0}$, which will
contradict (\ref{form55}).

{}From (\ref{0ik}) we have
\begin{equation} \label{form11}
a = b_1x_{\alpha}(-1)+c_1
\end{equation}
with
\begin{equation} \label{co1}
b_1 \in U(\bar{\goth{n}}_{-}) \; \; \mbox{and} \; \; c_1 \in
I_{\Lambda_{0}}.
\end{equation}
By Remark \ref{L(0)} we may and do assume that $b_1$ and $c_1$ are
doubly homogeneous; then $\mbox{wt} \; b_1= \mbox{wt} \; a-1$,
$\mbox{wt} \; c_1= \mbox{wt} \; a$, the charge of $b_1$ is one less
than that of $a$, and $c_1$ and $a$ have the same charge.

We now claim that
\begin{equation} \label{form12}
b_1 x_{\alpha}(-1) \in
I_{\Lambda_{0}}.
\end{equation}
Assume then that
\begin{equation} \label{formm}
b_1 x_\alpha(-1) \notin I_{\Lambda_{0}}.
\end{equation}
Then
\begin{equation} \label{b1notinspace}
b_1 \notin U(\bar{\goth{n}}_-)x_{\alpha}(-1);
\end{equation}
otherwise, $b_1x_{\alpha}(-1) \in
U(\bar{\goth{n}}_{-})x_{\alpha}(-1)^2 \subset I_{\Lambda_{0}}$.
By (\ref{deco}) we have a unique decomposition
\begin{equation} \label{form13}
b_1=r_1 x_{\alpha}(-1)+s_1, \ \
r_1 \in U(\bar{\goth{n}}_{-}), \ \ s_1 \in U(\bar{\goth{n}}_{\leq -2}),
\end{equation}
and $r_1$ and $s_1$ are doubly homogeneous, with $\mbox{wt} \;
r_1=\mbox{wt} \; b_1-1$, $\mbox{wt}\; s_1=\mbox{wt}\; b_1$, and
similarly for charge.  We have $s_1 \neq 0$ by (\ref{b1notinspace}).
We will use the vector $s_1$ to produce a contradiction.  We have
\begin{equation} \label{form14}
\tau^{-1}(b_1)=\tau^{-1}(r_1)x_{\alpha}(0)+\tau^{-1}(s_1) \; \;
\mbox{and} \; \; \tau^{-1}(s_1) \in U(\bar{\goth{n}}_{-}).
\end{equation}
Since
$b_1x_{\alpha}(-1) \cdot v_{\Lambda_{0}}^N= a \cdot
v_{\Lambda_{0}}^N-c_1 \cdot v_{\Lambda_{0}}^N \in \mbox{Ker} \;
\pi_{\Lambda_{0}}$,
$$
b_1 x_{\alpha}(-1) \cdot v_{\Lambda_{0}}=0,
$$
and so by (\ref{t1}),
$$
\tau^{-1} (b_1) \cdot v_{\Lambda_{1}} = 0.
$$
Thus (\ref{form14}) gives
$$
\tau^{-1} (s_1) \cdot v_{\Lambda_{1}} = 0,
$$
so that
\begin{equation} \label{form16}
\tau^{-1}(s_1) \cdot v_{\Lambda_{1}}^N \in \mbox{Ker} \;
\pi_{\Lambda_{1}}.
\end{equation}
By combining (\ref{t2}), (\ref{formm}) and (\ref{form13}) we also have
\begin{equation}
\widehat{e^{\alpha/2}_{(1, 1)}} (\tau^{-1} (s_1) \cdot
v_{\Lambda_{1}}^N)= s_1 x_{\alpha}(-1) \cdot
v_{\Lambda_{0}}^N=b_1x_{\alpha}(-1) \cdot
v_{\Lambda_0}^N-r_1x_{\alpha}(-1)^2 \cdot v_{\Lambda_0}^N \notin
I_{\Lambda_{0}} \cdot v_{\Lambda_{0}}^N,
\end{equation}
which by Lemma \ref{lemma2} implies
\begin{equation} \label{c2}
\tau^{-1} (s_1) \cdot v_{\Lambda_{1}}^N \notin I_{\Lambda_{1}}
\cdot v_{\Lambda_{1}}^N.
\end{equation}
Since $s_1$ is doubly homogeneous, so is $\tau^{-1} (s_1)$, and
\begin{equation} \label{w-formula}
\mbox{wt} \; \tau^{-1} (s_1) \leq \mbox{wt} \; s_1 = \mbox{wt} \; b_1
< \mbox{wt} \; a = L
\end{equation}
(note that if $s_1$ has charge 0, that is, is a constant, then
$\mbox{wt} \; \tau^{-1} (s_1) = \mbox{wt} \; s_1$).  Now
(\ref{form16}) and (\ref{c2}) together with the fact that $\tau^{-1}
(s_1)$ is doubly homogeneous of weight less than $L$ give us a
contradiction.  This proves our claim (\ref{form12}), and hence that
\begin{equation}
a= b_1x_{\alpha}(-1)+c_1 \in I_{\Lambda_0},
\end{equation}
which contradicts (\ref{form55}).  We have proved that $i$ cannot be
$0$ and we have thus established (\ref{sufftoshow}), completing the
proof of Theorem \ref{th2} for $k=1$. $\Box$

\begin{remark} \label{reasons}
\rm We have just proved Theorem \ref{th2} (formula (\ref{eqone})) for
$k=1$, by contradiction, in such a way that the assertion to be
contradicted, namely, (\ref{form1}), involves {\it both}
$W^N(\Lambda_0)$ and $W^N(\Lambda_1)$.  A different proof of this
theorem was given in \cite{CalLM1} (see the proof of Theorem 2.2),
where our argument proved the result for $W^N(\Lambda_0)$ and used
this result to prove the result for $W^N(\Lambda_1)$.  Also, the proof
given here does not use the space $W^N(\Lambda_1)'$ and related
``primed'' spaces, which played a crucial role in the proof of the
corresponding result in \cite{CalLM1}.  We have, however, included
information about such ``primed'' spaces in the present paper,
including conclusions about them in Theorems \ref{th1} and \ref{th2},
partly for reasons of comparison with the arguments in \cite{CalLM1}.
Our new argument for proving the $k=1$ case of Theorem \ref{th2}
naturally generalizes to $k \geq 1$ (see the proof below, which, while
it certainly reduces to the proof above when $k=1$, appears more
complicated in the greater generality), and it will also be
generalized in a different direction in subsequent work \cite{CalLM2}.
\end{remark}

We now generalize the $k=1$ proof to all $k \geq 1$.

{\it Proof of Theorem \ref{th2}:}
In view of (\ref{first inclusions}) it is sufficient to prove that
\begin{equation}\label{k_i}
\mbox{Ker} \; \pi_{\Lambda_{k,i}} \subset I_{\Lambda_{k,i}} \cdot
v_{\Lambda_{k,i}}^N \; \; \mbox{for all} \; \; i=0,\dots,k.
\end{equation}

Again we will prove this by contradiction.  Suppose then that there
exists $a \in U(\bar{\goth{n}}_{-})$ such that
\begin{equation} \label{contra}
a \cdot v_{\Lambda_{k,i}}^N \in \mbox{Ker} \; \pi_{\Lambda_{k,i}}
\; \; \; \mbox{but} \; \; \; a \cdot v_{\Lambda_{k,i}}^N \notin
I_{\Lambda_{k,i}} \cdot v_{\Lambda_{k,i}}^N \ \ {\rm for \ some}
\; \; i=0,\dots,k.
\end{equation}
By Remarks \ref{L(0)} and \ref{L(0)k} we may and do assume that $a$ is
doubly homogeneous.  Since $a$ is nonzero and in fact nonconstant (as
above), it has positive weight and charge.  Let
\begin{equation}\label{def-Lk}
L={\rm min} \{ {\rm wt} \; d \; | \; d \in U(\bar{\goth{n}}_-) \;
{\rm{doubly}}\;{\rm{ homogeneous}}\; {\rm{such \ that}} \
(\ref{contra}) \ {\rm{holds \ for}}\; d \} \; \; \; (>0).
\end{equation}
We further assume that $\mbox{wt} \; a =L$. Note that $i$ might be any
one or more of the indices from $0$ to $k$.  We shall show first that in
fact $i$ cannot be $k$.

Formulas (\ref{form2})--(\ref{wttau-1v0}) hold, exactly as in the
$k=1$ case.

Suppose that $i=k$, that is,
\begin{equation} \label{contraik}
a \cdot v_{\Lambda_{k,k}}^N \in \mbox{Ker} \; \pi_{\Lambda_{k,k}} \;
\; \; \mbox{but} \; \; \; a \cdot v_{\Lambda_{k,k}}^N \notin
I_{\Lambda_{k,k}} \cdot v_{\Lambda_{k,k}}^N \; \; \; (\mbox{that is,}
\; \; \; a \notin I_{\Lambda_{k,k}}),
\end{equation}
where $a$ is doubly homogeneous and $\mbox{wt} \; a =L$ (recall
(\ref{def-Lk})).  We will show that $\tau^{-1} (s_0)$ (recall
(\ref{form2})) is a doubly homogeneous element of
$U(\bar{\goth{n}}_-)$ whose weight is less than $L$ and which
satisfies (\ref{contra}).  We see that $s_0 \neq 0$, since $a \notin
U(\bar{\goth{n}}_{-})x_{\alpha}(-1)$, by (\ref{0ik}) and
(\ref{contraik}), and we know that $\mbox{wt} \; \tau^{-1} (s_0) <
\mbox{wt} \; a$.  {}From (\ref{form2}) and (\ref{contraik}) we obtain
$s_0 \cdot v_{\Lambda_{k,k}}^N \in \mbox{Ker} \; \pi_{\Lambda_{k,k}}$,
which is equivalent to $s_0 \cdot v_{\Lambda_{k,k}}=0$. Since
$$e^{\alpha/2}_{(k,0)}(\tau^{-1}(s_0) \cdot v_{\Lambda_{k,0}})= s_0
\cdot v_{\Lambda_{k,k}}= 0$$
(from (\ref{t1})) and since $e^{\alpha/2}_{(k,0)}$ is injective we
obtain
\be
\label{ikker}
\tau^{-1}(s_0) \cdot v_{\Lambda_{k,0}}^N \in {\rm
Ker} \; \pi_{\Lambda_{k,0}}.
\ee
Just as in the proof of the case $k=1$ we show that
\begin{equation} \label{ikI}
\tau^{-1} (s_0) \cdot v_{\Lambda_{k,0}}^N \notin I_{\Lambda_{k,0}}
\cdot v_{\Lambda_{k,0}}^N,
\end{equation}
and we have constructed a doubly homogeneous element $\tau^{-1}(s_0)$
of $U(\bar{\goth{n}}_-)$ satisfying (\ref{contra}) whose weight is
less than $\mbox{wt} \; a=L$.  This is a contradiction, and so $i$
cannot be $k$.

Now we may and do assume that
\begin{equation} \label{level-ki}
a \cdot v_{\Lambda_{k,i}}^N \in \mbox{Ker} \; \pi_{\Lambda_{k,i}} \;
\; \mbox{but} \; \; a \cdot v_{\Lambda_{k,i}}^N \notin
I_{\Lambda_{k,i}} \cdot v_{\Lambda_{k,i}}^N \; \; \mbox{for some} \;
\; i=0, \dots, k-1,
\end{equation}
where $a$ is doubly homogeneous of weight $L$ (recall
(\ref{def-Lk})). We now fix any one of the indices $i$
for which (\ref{level-ki}) holds.  Our next goal is to
show that $i$ cannot be $k-1$.

Since $a \cdot v_{\Lambda_{k,i}} \in \mbox{Ker} \; \pi_{\Lambda_{k,i}}$
we have
$$
a \cdot v_{\Lambda_{k,i}}=0 \; \; \mbox{in} \; \; W(\ki),
$$
and thus by Lemma \ref{kerns} we obtain
\begin{equation} \label{vk}
a \cdot v_{\Lambda_{k,k}}=0 \; \; \mbox{in} \; \; W(\Lambda_{k,k}),
\end{equation}
that is, $a \cdot v_{\Lambda_{k,k}}^N \in \mbox{Ker} \;
\pi_{\Lambda_{k,k}}$. The case we just proved (that $i$
cannot be $k$) thus gives us
$$
a \cdot v_{\Lambda_{k,k}}^N \in I_{\Lambda_{k,k}} \cdot
v_{\Lambda_{k,k}}^N,
$$
and so
$$a \in I_{\Lambda_{k,k}}.$$

Just as in (\ref{form11})--(\ref{co1}), we use (\ref{0ik}) to write
\begin{equation} \label{a-1}
a = b_1x_{\alpha}(-1)+c_1
\end{equation}
with
\begin{equation} \label{coo1}
b_1 \in U(\bar{\goth{n}}_{-}) \; \; \mbox{and} \; \; c_1 \in
I_{\Lambda_{k,0}},
\end{equation}
and by Remark \ref{L(0)} we may and do assume that $b_1$ and $c_1$ are
doubly homogeneous.  Then in fact $\mbox{wt} \;b_1=\mbox{wt} \; a -1$,
$\mbox{wt} \; c_1=\mbox{wt} \; a$, the charge of $b_1$ is one less
than that of $a$, and $c_1$ and $a$ have the same charge.

We next claim that
\begin{equation} \label{claim2}
b_1 x_{\alpha}(-1) \in
I_{\Lambda_{k,k-1}}.
\end{equation}
Suppose instead that
\begin{equation} \label{con}
b_1 x_\alpha(-1) \notin I_{\Lambda_{k,k-1}}.
\end{equation}
Then
\begin{equation} \label{nonclaim2}
b_1 \notin U(\bar{\goth{n}}_-)x_{\alpha}(-1)
\end{equation}
(otherwise,
$b_1x_{\alpha}(-1) \in U(\bar{\goth{n}}_{-})x_{\alpha}(-1)^2 \subset
I_{\Lambda_{k,k-1}}$).  We have a unique
decomposition
\be
\label{decompa1} b_1=r_1 x_{\alpha}(-1)+s_1, \ \ r_1
\in U(\bar{\goth{n}}_{-}), \ \ s_1 \in
U(\bar{\goth{n}}_{\leq -2})
\ee
by (\ref{deco}), and $r_1$ and $s_1$ are doubly homogeneous, with
$\mbox{wt} \; r_1=\mbox{wt} \; b_1-1$, $\mbox{wt} \; s_1=\mbox{wt}\;
b_1$, and similarly for charge.  Note that by (\ref{nonclaim2}) we
have $s_1 \neq 0$.  We also have
\begin{equation} \label{decompa1-tau}
\tau^{-1}(b_1)=\tau^{-1}(r_1)x_{\alpha}(0)+ \tau^{-1}(s_1) \; \;
\mbox{and} \; \; \tau^{-1} (s_1) \in U(\bar{\goth{n}}_{-}).
\end{equation}
Remark \ref{inclusions} and (\ref{first inclusions}) yield the inclusions
\begin{equation}\label{inclusionsformula}
I_{\Lambda_{k,0}} \cdot v_{\Lambda_{k,i}}^N \subset I_{\Lambda_{k,i}}
\cdot v_{\Lambda_{k,i}}^N \subset \mbox{Ker} \; \pi_{\Lambda_{k,i}},
\end{equation}
so that
$$b_1x_{\alpha}(-1) \cdot v_{\Lambda_{k,i}}^N =(a-c_1) \cdot
v_{\Lambda_{k,i}}^N \in \mbox{Ker} \; \pi_{\Lambda_{k,i}}$$
(recall (\ref{level-ki}), (\ref{a-1}) and (\ref{coo1})), and this is
equivalent to
$$
b_1 x_{\alpha}(-1) \cdot v_{\Lambda_{k,i}}=0.
$$
Now by Lemma \ref{kerns}
we obtain
\begin{equation} \label{b_1}
b_1x_{\alpha}(-1) \cdot v_{\Lambda_{k,k-1}}=0,
\end{equation}
and so (\ref{t1}) yields
$$
\tau^{-1}(b_1) \cdot v_{\Lambda_{k,1}} = 0.
$$
Hence from (\ref{decompa1-tau}) we get
\begin{equation} \label{b_1-ker}
\tau^{-1} (s_1) \cdot v_{\Lambda_{k,1}}^N =\tau^{-1}(b_1) \cdot
v_{\Lambda_{k,1}}^N \in \mbox{Ker} \; \pi_{\Lambda_{k,1}}.
\end{equation}
On the other hand, by (\ref{t2}) (using the fact that $\tau^{-1} (s_1)
\in U(\bar{\goth{n}}_{-})$), (\ref{con}) and (\ref{decompa1}) we also have
$$
\widehat{e^{\alpha/2}_{(k,1)}} (\tau^{-1} (s_1) \cdot
v_{\Lambda_{k,1}}^N)= s_1 x_{\alpha}(-1) \cdot v_{\Lambda_{k,k-1}}^N=
b_1x_{\alpha}(-1) \cdot v_{\Lambda_{k,k-1}}^N-r_1x_{\alpha}(-1)^2
\cdot v_{\Lambda_{k,k-1}}^N \notin I_{\Lambda_{k,k-1}} \cdot
v_{\Lambda_{k,k-1}}^N.
$$
Thus by using Lemma \ref{lemma2} we obtain
\begin{equation} \label{b_1-I}
\tau^{-1} (s_1 ) \cdot v_{\Lambda_{k,1}}^N \notin I_{\Lambda_{k,1}}
\cdot v_{\Lambda_{k,1}}^N.
\end{equation}
Just as in the proof of the case $k=1$ (recall (\ref{w-formula})) we
see that $\tau^{-1}(s_1)$ is doubly homogeneous of weight less than
$\mbox{wt} \; a=L$. We have obtained a contradiction by constructing
the doubly homogeneous element $\tau^{-1} (s_1)$ satisfying
(\ref{b_1-ker}) and (\ref{b_1-I}), and hence (\ref{contra}), whose
weight is less than $L$.  (Note that if $k=1$, we are {\it not}
claiming that $\tau^{-1} (s_1)$ also satisfies (\ref{level-ki}).)
This proves our claim (\ref{claim2}).

Thus
\begin{equation}
a= b_1x_{\alpha}(-1)+c_1 \in I_{\Lambda_{k,k-1}}
\end{equation}
with $b_1 \in U(\bar{\goth{n}}_{-})$ and $c_1 \in I_{\Lambda_{k,0}}$,
by (\ref{a-1}) and (\ref{eq0}), and we have shown that the index $i$
in (\ref{contra}) and in (\ref{level-ki}) cannot be $k-1$.  In
particular, if $k=1$ we are done.

Suppose then that $k \geq 2$.  Then we may and
do choose the index $i$ in (\ref{level-ki}) so that $0 \leq i \leq
k-2$.  We shall next show that this index $i$ cannot be $k-2$.  This
argument will be similar to the previous one, and it will make the
general pattern clear.

Since
\begin{equation} \label{ak-2}
a \cdot v_{\Lambda_{k,i}} =0 \; \; \mbox{for some} \; \; i=0,
\dots, k-2,
\end{equation}
by Lemma \ref{kerns} we get
\begin{equation}
a \cdot v_{\Lambda_{k,k-1}} =0,
\end{equation}
that is, $a \cdot v_{\Lambda_{k,k-1}}^N \in \mbox{Ker} \;
\pi_{\Lambda_{k,k-1}}$.  By the previous case (that $i$ cannot be
$k-1$),
\begin{equation}
a \cdot v_{\Lambda_{k,k-1}}^N \in I_{\Lambda_{k,k-1}} \cdot
v_{\Lambda_{k,k-1}}^N,
\end{equation}
so that
\begin{equation} \label{aIk-1}
a \in I_{\Lambda_{k,k-1}}.
\end{equation}
Thus from (\ref{0ik}) we obtain
\begin{equation} \label{a-2}
a=b_2x_{\alpha}(-1)^2+c_2 \; \; \mbox{with} \; \; b_2 \in
U(\bar{\goth{n}}_{-}) \; \; \mbox{and} \; \; c_2\in I_{\Lambda_{k,0}},
\end{equation}
and as usual, we may and do assume that $b_2$ and $c_2$ are doubly
homogeneous (by Remark \ref{L(0)}), so that $\mbox{wt} \;
b_2=\mbox{wt} \; a-2$, $\mbox{wt} \; c_2= \mbox{wt} \; a$, the charge
of $b_2$ is two less that of $a$, and $c_2$ and $a$ have the same
charge.

We now prove by contradiction that
\begin{equation} \label{claim3}
b_2x_{\alpha}(-1)^2 \in I_{\Lambda_{k, k-2}}
\end{equation}
(cf. (\ref{claim2})): If instead
\begin{equation} \label{notclaim3}
b_2x_{\alpha}(-1)^2 \notin I_{\Lambda_{k, k-2}},
\end{equation}
then $b_2 \notin U(\bar{\goth{n}}_{-})x_{\alpha}(-1)$
(cf. (\ref{nonclaim2})), and thus we have a unique decomposition
\begin{equation} \label{b2r2s2}
b_2=r_2x_{\alpha}(-1)+s_2, \; \; r_2 \in U(\bar{\goth{n}}_{-}), \;
\; 0 \neq s_2 \in U(\bar{\goth{n}}_{\leq -2}),
\end{equation}
and $r_2$ and $s_2$ are doubly homogeneous, with $\mbox{wt} \; r_2=
\mbox{wt} \; b_2-1$, $\mbox{wt} \; s_2= \mbox{wt} \; b_2$, and
similarly for charge (as in (\ref{decompa1})).  We follow the argument
of (\ref{decompa1-tau})--(\ref{b_1-I}): We apply $\tau^{-1}$ to
(\ref{b2r2s2}).  Since (\ref{inclusionsformula}) still holds, we
obtain that
$$b_2x_{\alpha}(-1)^2 \cdot v_{\Lambda_{k,i}} =(a-c_2) \cdot
v_{\Lambda_{k,i}} = 0,$$
which gives
$$b_2x_{\alpha}(-1)^2 \cdot v_{\Lambda_{k,k-2}} = 0$$
by Lemma \ref{kerns}.  Thus by (\ref{t1}) we get
$$
\tau^{-1}(b_2) \cdot v_{\Lambda_{k,2}} = 0,
$$
and so
\begin{equation} \label{b_2-ker}
\tau^{-1} (s_2) \cdot v_{\Lambda_{k,2}}^N =\tau^{-1}(b_2) \cdot
v_{\Lambda_{k,2}}^N \in \mbox{Ker} \; \pi_{\Lambda_{k,2}}.
\end{equation}
Using (\ref{t2}), the fact that $\tau^{-1} (s_2) \in
U(\bar{\goth{n}}_{-})$, (\ref{notclaim3}) and (\ref{b2r2s2}), we also
obtain
$$
(2!) \, \widehat{e^{\alpha/2}_{(k,2)}} (\tau^{-1} (s_2) \cdot
v_{\Lambda_{k,2}}^N)= b_2x_{\alpha}(-1)^2 \cdot
v_{\Lambda_{k,k-2}}^N-r_2x_{\alpha}(-1)^3 \cdot v_{\Lambda_{k,k-2}}^N
\notin I_{\Lambda_{k,k-2}} \cdot v_{\Lambda_{k,k-2}}^N,
$$
and so by Lemma \ref{lemma2},
\begin{equation} \label{b_2-I}
\tau^{-1} (s_2 ) \cdot v_{\Lambda_{k,2}}^N \notin I_{\Lambda_{k,2}}
\cdot v_{\Lambda_{k,2}}^N.
\end{equation}
Just as in the proof above, $\tau^{-1}(s_2)$ is a doubly homogeneous
element satisfying (\ref{b_2-ker}) and (\ref{b_2-I}) and hence
(\ref{contra}) (but not necessarily (\ref{level-ki})) and of weight
less than $L$.  This proves (\ref{claim3}).

By (\ref{a-2}), (\ref{claim3}) and (\ref{eq0}) we now have
\begin{equation}
a=b_2x_{\alpha}(-1)^2+c_2 \in I_{\Lambda_{k,k-2}}
\end{equation}
with $b_2 \in U(\bar{\goth{n}}_{-})$ and $c_2\in I_{\Lambda_{k,0}}$,
and this proves that $i$ cannot be $k-2$.  In particular, we are done
if $k=2$.

Now we give the general inductive step.  Fix $m \geq 1$ and assume
that the assertion of Theorem \ref{th2} has been proved for $k =
1,2,\dots, m$ and that $i$ in (\ref{contra}) (or in (\ref{level-ki}))
cannot be $k, k-1, \dots, k-m$.  We shall show that if $k \geq m+1$,
then the index $i$ cannot be $k-(m+1)$ either, and that in particular,
the assertion of Theorem \ref{th2} thus holds for $k=m+1$.  This will
complete the proof of the theorem.

Suppose then that $k \geq m+1$ and that the index $i$ in
(\ref{level-ki}) is such that $0 \leq i \leq k-(m+1)$.  To show that
this index $i$ in fact cannot be $k-(m+1)$, we first observe that
exactly as in (\ref{ak-2})--(\ref{aIk-1}) we have
$$
a \in I_{\Lambda_{k,k-m}},
$$
and so from (\ref{0ik}) we see that
\begin{equation} \label{a-m+1}
a=b_{m+1}x_{\alpha}(-1)^{m+1}+c_{m+1} \; \; \mbox{with} \; \; b_{m+1} \in
U(\bar{\goth{n}}_{-}) \; \; \mbox{and} \; \; c_{m+1}\in I_{\Lambda_{k,0}}.
\end{equation}
Again, as above, we may and do assume that $b_{m+1}$ and $c_{m+1}$ are
doubly homogeneous (by Remark \ref{L(0)}); then $\mbox{wt} \;
b_{m+1}=\mbox{wt} \; a-(m+1)$, $\mbox{wt} \; c_{m+1}= \mbox{wt} \; a$,
the charge of $b_{m+1}$ is $m+1$ less that of $a$, and $c_{m+1}$ and
$a$ have the same charge.

Exactly as in (\ref{claim3})--(\ref{b_2-I}), we obtain by
contradiction that
\begin{equation}\label{claim4}
b_{m+1}x_{\alpha}(-1)^{m+1} \in I_{\Lambda_{k, k-(m+1)}}:
\end{equation}
Assume that
$$
b_{m+1}x_{\alpha}(-1)^{m+1} \notin I_{\Lambda_{k, k-(m+1)}}.
$$
In place of formula (\ref{b2r2s2}), we now have the unique
decomposition
$$
b_{m+1}=r_{m+1}x_{\alpha}(-1)+s_{m+1}, \; \; r_{m+1} \in U(\bar{\goth{n}}_{-}), \;
\; 0 \neq s_{m+1} \in U(\bar{\goth{n}}_{\leq -2}),
$$
with $r_{m+1}$ and $s_{m+1}$ doubly homogeneous, $\mbox{wt} \;
r_{m+1}= \mbox{wt} \; b_{m+1}-1$, $\mbox{wt} \; s_{m+1}= \mbox{wt} \;
b_{m+1}$, and similarly for charge.  As in formula
(\ref{decompa1-tau}) we now have
$$
\tau^{-1}(b_{m+1})=\tau^{-1}(r_{m+1})x_{\alpha}(0)+ \tau^{-1}(s_{m+1})
\; \; \mbox{and} \; \; \tau^{-1} (s_{m+1}) \in U(\bar{\goth{n}}_{-}).
$$
By (\ref{inclusionsformula}) we obtain
$$
b_{m+1}x_{\alpha}(-1)^{m+1} \cdot v_{\Lambda_{k,i}} =(a-c_{m+1}) \cdot
v_{\Lambda_{k,i}} = 0,
$$
so that
$$
b_{m+1}x_{\alpha}(-1)^{m+1} \cdot v_{\Lambda_{k,k-({m+1})}} = 0,
$$
by Lemma \ref{kerns}, and so (\ref{t1}) gives
$$
\tau^{-1}(b_{m+1}) \cdot v_{\Lambda_{k,{m+1}}} = 0.
$$
Thus
\begin{equation}\label{b_m+1-ker}
\tau^{-1} (s_{m+1}) \cdot v_{\Lambda_{k,{m+1}}}^N =\tau^{-1}(b_{m+1}) \cdot
v_{\Lambda_{k,{m+1}}}^N \in \mbox{Ker} \; \pi_{\Lambda_{k,{m+1}}}.
\end{equation}
Since $\tau^{-1} (s_{m+1}) \in U(\bar{\goth{n}}_{-})$, we can use
(\ref{t2}), and exactly as above we find that
$$
(m+1)! \, \widehat{e^{\alpha/2}_{(k,{m+1})}} (\tau^{-1} (s_{m+1}) \cdot
v_{\Lambda_{k,{m+1}}}^N)= b_{m+1}x_{\alpha}(-1)^{m+1} \cdot
v_{\Lambda_{k,k-(m+1)}}^N-r_{m+1}x_{\alpha}(-1)^{m+2} \cdot
v_{\Lambda_{k,k-(m+1)}}^N,
$$
so that
$$
\widehat{e^{\alpha/2}_{(k,{m+1})}} (\tau^{-1} (s_{m+1}) \cdot
v_{\Lambda_{k,{m+1}}}^N) \notin I_{\Lambda_{k,k-(m+1)}} \cdot
v_{\Lambda_{k,k-(m+1)}}^N.
$$
Thus by Lemma \ref{lemma2},
\begin{equation}\label{b_m+1-I}
\tau^{-1} (s_{m+1} ) \cdot v_{\Lambda_{k,{m+1}}}^N \notin I_{\Lambda_{k,{m+1}}}
\cdot v_{\Lambda_{k,{m+1}}}^N.
\end{equation}
Since $\tau^{-1}(s_{m+1})$ is a doubly homogeneous element satisfying
(\ref{b_m+1-ker}) and (\ref{b_m+1-I}) and thus (\ref{contra}) (but not
necessarily (\ref{level-ki})) and of weight less than $L$, we have
proved (\ref{claim4}).

Hence from (\ref{a-m+1}), (\ref{claim4}) and (\ref{eq0}) we finally
obtain
\begin{equation}
a=b_{m+1}x_{\alpha}(-1)^{m+1}+c_{m+1} \in I_{\Lambda_{k,k-(m+1)}},
\end{equation}
proving that $i$ cannot be $k-(m+1)$ and thus proving Theorem
\ref{th2}.  $\; \; \; \Box$

\begin{remark}
\rm The first part of the proof, in which we showed that $i$ cannot be
$k$, is actually essentially the same argument as the successive
arguments showing that $i$ cannot be $k-1, k-2$, and so on.
\end{remark}

\begin{remark}
\rm As an immediate consequence of Theorem \ref{th1}, we see that any
nonzero doubly homogeneous element $a \in U(\bar{\goth{n}}_{-})$ such
that $a \in \mbox{Ker} \; f_{\Lambda_{k,0}}= I_{\Lambda_{k,0}} $ has
charge at least $k+1$; that is, no nonzero linear combination of
monomials $x_{\alpha}(m_1) \cdots x_{\alpha}(m_r)$ with $r \leq k$ and
each $m_i<0$ belongs to $\mbox{Ker} \; f_{\Lambda_{k,0}}$.  We observe
similarly that any homogeneous element of charge $k+1$ that lies in
$\mbox{Ker} \; f_{\Lambda_0}$ is a multiple of $R_t^0$ for some $t
\geq k+1$.
\end{remark}

\section{Another reformulation}
\setcounter{equation}{0}

Generalizing the last section of \cite{CalLM1}, we shall finally give
a further reformulation of the $i=0$ case of Theorem \ref{th2},
formula (\ref{eqone}), in terms of principal ideals of vertex
(operator) algebras.  As in \cite{CalLM1}, we shall invoke \cite{LL}
for material on ideals of vertex (operator) algebras and on vertex
operator algebra and module structure on generalized Verma modules.

The generalized Verma module $N(\Lambda_{k,0})$ has a natural
structure of vertex operator algebra, with vertex operator map
\begin{eqnarray} \nonumber
Y(\cdot, x): N(\Lambda_{k,0}) & \longrightarrow & \mbox{End} \;
N(\Lambda_{k,0}) [[x, x^{-1}]] \nonumber \\
v & \mapsto & Y(v, x)= \sum_{m \in \mathbb{Z}} v_{m} x^{-m-1}
\nonumber
\end{eqnarray}
satisfying the conditions given in Theorem 6.2.18 of \cite{LL}, with
$v_{\Lambda_{k,0}}^N$ as vacuum vector.  The conformal vector gives
rise to the Virasoro algebra operators $L(m)$, $m \in \mathbb{Z}$,
including the operator $L(0)$ used above.  Also, $N(\ki)$ for $0 \leq
i \leq k$ is naturally a module for the vertex operator algebra
$N(\Lambda_{k,0})$, as described in Theorem 6.2.21 of \cite{LL}.

Just as in \cite{CalLM1}, $W^N(\Lambda_{k,0})$ is a vertex subalgebra
of $N(\Lambda_{k,0})$ and $W^N(\ki)$ is a
$W^N(\Lambda_{k,0})$-submodule of $N(\ki)$ for $0 \leq i \leq k$.
Also, $L(0)$ preserves $W^N(\ki)$ for $0 \leq i \leq k$ and $L(-1)$
preserves only $W^N(\Lambda_{k,0})$.

We recall from Section 3 the natural surjective
$\widehat{\goth{g}}$-module maps
\begin{eqnarray} \label{Pi}
\Pi_{\Lambda_{k,i}} : N(\ki) & \longrightarrow & L(\ki) \\
a \cdot v_{\Lambda_{k,i}}^N & \mapsto & a \cdot v_{\Lambda_{k,i}}, \;
\; \; \; a \in U(\widehat{\goth{g}})\nonumber
\end{eqnarray}
and their kernels
\begin{equation} \label{ker}
N^1(\ki)= \mbox{Ker} \; \Pi_{\Lambda_{k,i}},
\end{equation}
for $0 \leq i \leq k$.  Then $N^1(\ki)$ is the unique maximal proper
($L(0)$-graded) $\widehat{\goth{g}}$-submodule of $N(\ki)$ and
$$
N^1(\Lambda_{k,i}) = U(\widehat{\goth{g}}) x_{\alpha}(-1)^{k-i+1}\cdot
v_{\Lambda_{k,i}}^N= U(\mathbb{C}x_{-\alpha} \oplus \goth{g} \otimes
t^{-1} \mathbb{C}[t^{-1}]) x_{\alpha}(-1)^{k-i+1} \cdot v_{\Lambda_{k,i}}^N
$$
for $ 0 \leq i \leq k$ (cf. \cite{K}, \cite{LL}).

As in \cite{CalLM1}, a {\it principal ideal} of a vertex (operator)
algebra is an ideal generated by a single element.  The following
result, which generalizes Proposition 4.1 in \cite{CalLM1} and which
is proved the same way, says that $N^1(\Lambda_{k,0})$ is the
principal ideal of $N(\Lambda_{k,0})$ generated by the ``null vector''
$x_{\alpha}(-1)^{k+1} \cdot v_{\Lambda_{k,0}}^N$:

\begin{proposition} \label{pro}
The space $N^1(\Lambda_{k,0})$ is the ideal of the vertex operator
algebra $N(\Lambda_{k,0})$ generated by $x_{\alpha}(-1)^{k+1}\cdot
v_{\Lambda_{k,0}}^N$. $\; \; \; \Box$
\end{proposition}

The kernels of the restrictions $\pi_{\Lambda_{k,i}}$ of the maps
(\ref{Pi}) to the principal subspaces $W^N(\ki)$ (recall (\ref{pi}))
are
\begin{equation} \label{newer}
\mbox{Ker} \; \pi_{\Lambda_{k,i}}= N^1(\ki) \cap
W^N(\ki)
\end{equation}
for $0 \leq i \leq k$. As in Remark 4.2 in \cite{CalLM1} we have that
$\mbox{Ker} \; \pi_{\Lambda_{k,0}}=N^1(\Lambda_{k,0}) \cap
W^N(\Lambda_{k,0})$, which equals $I_{\Lambda_{k,0}} \cdot
v_{\Lambda_{k,0}}^N$ by Theorem \ref{th2}, is an ideal of the vertex
algebra $W^N( \Lambda_{k,0})$.  Moreover, generalizing Proposition 4.2
of \cite{CalLM1} and using essentially the same proof, we have that
this ideal is also a principal ideal, generated by the same null
vector:

\begin{proposition} \label{4}
The space $I_{\Lambda_{k,0}} \cdot v_{\Lambda_{k,0}}^N$ is the ideal
of the vertex algebra $W^N(\Lambda_{k,0})$ generated by
$x_{\alpha}(-1)^{k+1} \cdot v_{\Lambda_{k,0}}^N$. $\; \; \; \Box$
\end{proposition}

Again as in \cite{CalLM1} we write $(v)_V$ for the ideal generated by
an element $v$ of a vertex (operator) algebra $V$.  Combining
Propositions \ref{pro} and \ref{4} with Theorem \ref{th2}, we have
obtained a reformulation of the $i=0$ case of Theorem \ref{th2},
formula (\ref{eqone}), genearalizing Theorem 4.1 of \cite{CalLM1}:

\begin{theorem} \label{th3}
For every $k > 0$,
\begin{equation}
{\rm Ker} \; \pi_{\Lambda_{k, 0}}=(x_{\alpha}(-1)^{k+1} \cdot
v_{\Lambda_{k,0}}^N ) _{N(\Lambda_{k,0})} \cap W^N( \Lambda_{k,0})=
(x_{\alpha}(-1)^{k+1} \cdot v_{\Lambda_{k,0}}^N)_{W^N(
\Lambda_{k,0})}.
\end{equation}
In particular, the intersection with the vertex subalgebra $W^N(
\Lambda_{k,0})$ of the principal ideal of $N( \Lambda_{k,0})$
generated by the null vector $x_{\alpha}(-1) ^{k+1} \cdot
v_{\Lambda_{k,0}}^N$ coincides with the principal ideal of the vertex
subalgebra $W^N( \Lambda_{k,0})$ generated by the same null
vector. $\; \; \; \Box$
\end{theorem}

\noindent {\small \sc Department of Mathematics, Rutgers University,
Piscataway, NJ 08854} \\
\noindent Current address: \\
\noindent{\small \sc Department of Mathematics,  Ohio State University,
Columbus, OH 43210} \\
{\em E--mail address}: calinescu@math.ohio-state.edu\\
On leave from the Institute of Mathematics of the Romanian Academy.\\
\vspace{.1in}

\noindent {\small \sc Department of Mathematics, Rutgers University,
Piscataway, NJ 08854} \\ {\em E--mail address}:
lepowsky@math.rutgers.edu \\
\vspace{.1in}

\noindent {\small \sc Department of Mathematics and Statistics,
University at Albany (SUNY), Albany, NY 12222} \\ {\em E--mail
address}: amilas@math.albany.edu

\end{document}